\documentclass[article]{elsarticle}

\usepackage{lineno,hyperref}
\usepackage{color}
\usepackage{caption}
\usepackage{subcaption}
\modulolinenumbers[1]
\usepackage{amsmath}

\usepackage{amsmath,amsfonts,bm}

\newcommand{\R}{\mathbb{R}}

\def\eqref#1{equation~\ref{#1}}

\def\1{\bm{1}}

\def\rva{{\mathbf{a}}}

\def\rvv{{\mathbf{v}}}

\def\rvz{{\mathbf{z}}}

\def\rmK{{\mathbf{K}}}

\DeclareMathAlphabet{\mathsfit}{\encodingdefault}{\sfdefault}{m}{sl}
\SetMathAlphabet{\mathsfit}{bold}{\encodingdefault}{\sfdefault}{bx}{n}

\def\sR{{\mathbb{R}}}

\DeclareMathOperator*{\argmin}{arg\,min}

\usepackage{booktabs}

\usepackage{tcolorbox} 
\usepackage{xcolor}
\usepackage{multicol}
\usepackage{algorithm}
\usepackage{algorithmic}
\usepackage{amssymb} 
\usepackage{amsmath} 
\usepackage{wrapfig}
\usepackage{graphicx}
\usepackage{multirow}
\usepackage{colortbl,booktabs} 
\usepackage{lipsum}   
\usepackage{amsthm} 
\usepackage{amsmath} 
\usepackage{arydshln} 
\usepackage{caption} 
\usepackage{comment}

 \newtheorem{theorem}{Theorem}[section]

 \newtheorem{remark}{Remark}[section]

\journal{Journal of \LaTeX\ Templates}

\begin{document}

\begin{frontmatter}

\title{Random Features for High-Dimensional Nonlocal Mean-Field Games}

\author[UCLACS]{Sudhanshu Agrawal\corref{cofiauthor}}
\author[UCLA]{Wonjun Lee\corref{cofiauthor}}
\author[Mines]{Samy Wu Fung\corref{cor}}\ead{swufung@mines.edu}
\author[UCLA]{Levon Nurbekyan}
\address[UCLACS]{Department of Computer Science, University of California, Los Angeles}
\address[UCLA]{Department of Mathematics, University of California, Los Angeles}
\address[Mines]{Department of Applied Mathematics and Statistics, Colorado School of Mines}
\cortext[cofiauthor]{Co-first Author}
\cortext[cor]{Corresponding Author}





\begin{abstract}
We propose an efficient solution approach for high-dimensional nonlocal mean-field game (MFG) systems based on the Monte Carlo approximation of interaction kernels via random features. We avoid costly space-discretizations of interaction terms in the state-space by passing to the feature-space. This approach allows for a seamless mean-field extension of virtually any single-agent trajectory optimization algorithm. Here, we extend the direct transcription approach in optimal control to the mean-field setting. We demonstrate the efficiency of our method by solving MFG problems in high-dimensional spaces which were previously out of reach for conventional non-deep-learning techniques.
\end{abstract}


\end{frontmatter}

\section{Introduction}


We propose a computational framework for solving mean-field game (MFG) systems of the form
\begin{equation}\label{eq:mfg_main}
    \begin{cases}
    -\partial\phi(t,x)+H(t,x,\nabla \phi(t,x))= \int_{\sR^d} K(x,y) d\rho(t,y) ~ \text{in} ~ (0,T)\times \sR^d,\\
    \partial_t \rho(t,x)-\nabla \cdot \left( \rho(t,x) \nabla_p H(t,x,\nabla \phi(t,x))\right)=0 ~ \text{in} ~ (0,T)\times \sR^d,\\
    \rho(0,x)=\rho_0(x),\quad \phi(T,x)=\psi(x) ~ \text{in}~ \sR^d,
    \end{cases}
\end{equation}
based on random features from kernel machines. The partial differential equation (PDE) above describes an equilibrium configuration of a noncooperative differential game with a continuum of agents. An individual agent faces a cost
\begin{equation}\label{eq:ind_cost}
    \phi(t,x)=\inf_{z(t)=x} \int_t^T \left\{L(s,z(s),\dot{z}(s))+\int_{\R^d} K(z(s),y) d\rho(s,y)\right\} ds+\psi(z(T)),
\end{equation}
where the \textit{Lagrangian (running cost)} $L$ and the \textit{Hamiltonian} $H$ are related by the \textit{Legendre transform}:
\begin{equation}\label{eq:L-H}
\begin{split}
    L(t,x,v)=&\sup_{p \in \sR^d}\{ -v\cdot p - H(t,x,p) \},\\
    H(t,x,p)=&\sup_{v \in \sR^d}\{-v\cdot p - L(t,x,v) \}. 
\end{split}
\end{equation}
Furthermore, $\rho(t,\cdot)$ represents the distribution of all agents in the state-space at time $t$, and the term
\begin{equation}\label{eq:interaction}
    f(x,\rho(t,\cdot))=\int_{\sR^d} K(x,y) d\rho(t,y)
\end{equation}
models the influence of the population on an individual agent. Finally, $\psi$ in \ref{eq:ind_cost} represents the terminal cost paid by agents at terminal time $T$, and $\rho_0$ is the initial distribution of the population. Note that \ref{eq:interaction} assumes a \textit{nonlocal interaction} of an individual agent with the population. If, for instance, we had
\begin{equation*}
    f(x,\rho(t,\cdot))=c~\rho(t,x)^\kappa \quad \text{or} \quad f(x,\rho(t,\cdot))=c~\log \rho(t,x)
\end{equation*}
then the interaction would be \textit{local}. In this paper, we only consider nonlocal interactions in \ref{eq:interaction}.

In an equilibrium, individual agents cannot unilaterally improve their costs based on their belief about the state-space distribution of the population. This Nash equilibrium principle leads to the Hamilton-Jacobi-Bellman (HJB) PDE in~\ref{eq:mfg_main}. Furthermore, the evolution of the state-space distribution of the population corresponding to their optimal actions must coincide with their belief about population distribution. This consistency principle leads to the continuity equation in \ref{eq:mfg_main}.

The MFG framework, introduced by M. Huang, P. Caines, R. Malham\'{e} \cite{HCM06,HCM07} and P.-L. Lions, J.-M. Lasry \cite{LasryLions06a,LasryLions06b,LasryLions2007}, is currently an active field with applications in economics~\cite{moll14,moll17,gueant2011mean,gomes2015economic}, finance~\cite{caines17,cardialiaguet2018,jaimungal19,moll17}, industrial engineering~\cite{paola19,kizikale19,gomes2018electricity}, swarm robotics~\cite{liu2018mean,elamvazhuthi2019mean,kang21jointsensing,kang21taskselection}, epidemic modelling~\cite{lee2020controlling,chang2020game} and data science~\cite{weinan2019mean,guo2019learning,carmona2019linear}. For comprehensive exposition of MFG systems we refer to \cite{LasryLions2007,CardaNotes,mfgCIME} for nonlocal couplings, \cite{gomes_book16,cirant19,mfgCIME} for local couplings,  \cite{delarue18a,delarue18b} for a probabilistic approach, \cite{bensoussan2013} for infinite-dimensional control approach, \cite{carda19,gangbo2021mean} for the master equation, and \cite{achdou:hal-03408825} for the control on the acceleration. For the mathematical analysis of \ref{eq:mfg_main} we refer to \cite{LasryLions2007,CardaNotes,mfgCIME}.

In this paper, we develop a computational method for \ref{eq:mfg_main} based on kernel expansion framework introduced in \cite{nurbekyan18,nursaude18,liu2020computational,liu2020splitting}. The key idea is to build an approximation
\begin{equation}\label{eq:kernel_approx}
    K(x,y) \approx K_r(x,y)=\sum_{i,j=1}^r k_{ij}\zeta_i(x) \zeta_j(y),
\end{equation}
where $\{\zeta_i\}_{i=1}^r$ and $(k_{ij})$ are suitably chosen basis functions and expansions coefficients, and consider an approximate system
\begin{equation}\label{eq:mfg_approx}
    \begin{cases}
    -\partial_t\phi(t,x)+H(t,x,\nabla \phi(t,x))= \int_{\sR^d} K_r(x,y) d\rho(t,y) ~ \text{in} ~ (0,T)\times \sR^d,\\
    \partial_t \rho(t,x)-\nabla \cdot \left( \rho(t,x) \nabla_p H(t,x,\nabla \phi(t,x))\right)=0 ~ \text{in} ~ (0,T)\times \sR^d,\\
    \rho(0,x)=\rho_0(x),\quad \phi(T,x)=\psi(x) ~ \text{in}~ \sR^d.
    \end{cases}
\end{equation}
The structure of $K_r$ allows for an efficient discretization of the interaction term $\int_{\R^d} K_r(x,y) d \rho(t,y)$ in the \textit{feature space}. Indeed, introducing \textit{unknown coefficients} $a(t)=(a_1(t),a_2(t),\cdots,a_r(t))$ we can rewrite \ref{eq:mfg_approx} as
\begin{equation}\label{eq:a_inclusion}
    0= \rmK^{-1} a(t)  - \frac{\delta}{\delta a(t)} \int_{\sR^d} \phi_a(0,x) d\rho_0(x), 
\end{equation}
where $\rmK=(k_{ij})_{i,j=1}^r$, and $\phi_a$ is the viscosity solution of
\begin{equation}\label{eq:phi_a}
    \begin{cases}
    -\partial_t\phi(t,x)+H(t,x,\nabla \phi(t,x))= \sum_{i=1}^r a_i(t) \zeta_i(x) ~ \text{in} ~ (0,T)\times \sR^d,\\
    \phi(T,x)=\psi(x) ~ \text{in}~ \sR^d.
    \end{cases}
\end{equation}
We provide a formal derivation of the equivalence between \ref{eq:mfg_approx} and \ref{eq:a_inclusion} in the Appendix and refer to \cite{nursaude18} for more details. When $\rmK$ is symmetric, \ref{eq:a_inclusion} reduces to an optimization problem
\begin{equation}\label{eq:a_optimization}
    \inf_a \frac{1}{2}\int_0^T \sum_{i,j=1}^r (\rmK^{-1})_{ij} a_i(t) a_j(t) dt  - \int_{\sR^d} \phi_a(0,x) d\rho_0(x). 
\end{equation}

The key advantage in our approach is that $(a_i)$ contain all information about the population interaction, and there is no need for a costly space discretization of $f$ in \ref{eq:interaction}. Indeed, the approximation in \ref{eq:kernel_approx} yields an approximation of the interaction \textit{operator}
\begin{equation*}
    \int_{\sR^d} K(x,y) d\rho(t,y) \approx \int_{\sR^d} K_r(x,y) d\rho(t,y)=\sum _{i=1}^r \zeta_i(x) \underbrace{\sum_{j=1}^r k_{ij}  \int_{\sR^d} \zeta_j(y) d\rho(t,y)}_{a_i(t)}
\end{equation*}
that is independent of the space-discretization. Moreover, for fixed $r$, the computational cost of calculating the approximate interaction term in space-time is $O(r^2N_t+2r N_x N_t)$, where $N_t$ is the time-discretization, and $N_x$ is the space-discretization or number of trajectories or agents in the Lagrangian setting. In contrast, direct calculation of the interaction term yields an $O(N_x^2 N_t)$ computational cost. This dimension reduction provides a significant computational gain when $r$ is moderate.

There is a complete flexibility in the choice of basis functions $\{\zeta_i\}$. In \cite{nursaude18}, the authors considered problems in periodic domains and used classical trigonometric polynomials. Furthermore, in \cite{liu2020computational,liu2020splitting} the authors drew connections with kernel methods in machine learning and used polynomial and quasi-polynomial features for $\{\zeta_i\}$.

Our key contribution is to build on the connection with kernel methods in machine learning and construct $\{\zeta_i\}$ using \textit{random features} \cite{RahimiR07}. The advantage of using random features for a suitable class of kernels is the simplicity and speed of the generation of basis functions, including in high dimensions. Moreover, $\rmK$ in \ref{eq:kernel_approx} reduces to an identity matrix which renders extremely simple update rules for $(a_i)$ in iterative solvers of \ref{eq:a_inclusion} and \ref{eq:a_optimization}.

We demonstrate the efficiency of our approach by solving crowd-motion-type MFG problems in up to $d=100$ dimensions. To the best of our knowledge, this is the first instance such high-dimensional MFG are solved without deep learning techniques. Our algorithm is inspired by the primal-dual algorithm in \cite{nursaude18}, except that here we use random features instead of trigonometric polynomials. The primal step consists of trajectory optimization, whereas the dual step updates nonlocal variables $(a_i)$. Modeling nonlocal interactions by $(a_i)$ \textit{decouples} primal updates for the agents, which would not be possible using a direct discretization of the interaction term. Hence, one can take advantage of parallelization techniques within primal updates. We refer to Section \ref{sec:algorithm} for more details.

For related work on numerical methods for nonlocal MFG we refer to~\cite{hadi17,hadi17b,hadi19,bonnans2021generalized} for game theoretic approach,~\cite{silva12,silva14,silva15,carlini18} for semi-Lagrangian schemes,~\cite{lin2021alternating} for deep learning approach, and~\cite{li2021multiscale} for a multiscale method. In all of these methods the nonlocal terms are discretized directly in the state-space. \textcolor{black}{Somewhat related work is~\cite{mou2021numerical} where the authors approximate the solutions in reproducing kernel Hilbert spaces (RKHS) and Fourier spaces. The critical difference is that we use features to approximate the interaction terms, not the solutions.} Finally, for a comprehensive exposition of numerical methods for other types of MFG systems we refer to~\cite{mfgCIME}.

The rest of the paper is organized as follows. In Section \ref{sec:coeffs} we present the kernel expansion framework. Next, in Section \ref{sec:rfeatures}, we show how to construct basis functions based on random features. Section \ref{sec:algorithm} contains the description of our algorithm. Finally, we present numerical results in Section \ref{sec:numerics}. We provide an implementation\footnote{code can be found in \url{https://github.com/SudhanshuAgrawal27/HighDimNonlocalMFG}} written in the Julia language~\cite{bezanson2017julia}.

\section{The method of coefficients}\label{sec:coeffs}

One can adapt the results in \cite{nursaude18} to the non-periodic setting relying on the analysis in \cite{CardaNotes} and prove the following theorem.
\begin{theorem}
A pair $(\phi,\rho)$ is a solution for the MFG \ref{eq:mfg_approx} if and only if there exist $a=(a_1,a_2,\cdots,a_r) \in C([0,T];\R^r)$ such that \ref{eq:a_inclusion} holds. Moreover, when $\rmK$ is symmetric, \ref{eq:a_inclusion} reduces to \ref{eq:a_optimization}. Finally, when $\rmK$ is positive-definite \ref{eq:a_optimization} is a convex program.
\end{theorem}

Next, we need a formula to calculate the gradient of the objective function in \ref{eq:a_optimization}. Again, adapting results in \cite{nursaude18} to the non-periodic setting one can prove the following theorem.

\begin{theorem}
The functional $a\mapsto \int_{\sR^d} \phi_a(0,x) d\rho_0(x)$ is convex and Fr\'{e}chet differentiable everywhere. Moreover,
\begin{equation}\label{eq:grad_intphi_a}
    \frac{\delta}{\delta a_i(t)} \int_{\sR^d} \phi_a(0,x) d\rho_0(x)=\int_{\sR^d} \zeta_i(z_{x,a}(t)) d\rho_0(x),
\end{equation}
where $z_{x,a}$ is an optimal trajectory for the optimal control problem
\begin{equation}\label{eq:z_ax}
    \phi_a(t,x)=\inf_{z(t)=x} \int_t^T \left\{L(s,z(s),\dot{z}(s))+
    \sum_{i=1}^r a_i(s)\zeta_i(z(s))\right\} ds+\psi(z(T)).
\end{equation}
\end{theorem}

We do not specify precise assumptions on the data in these previous theorems and refer to \cite{nursaude18,CardaNotes} for more details since the theoretical analysis of \ref{eq:mfg_main} and \ref{eq:mfg_approx} is out of the scope of the current paper. Nevertheless, these theorems are valid for typical choices such as
\begin{equation*}
    L(t,x,v)=\frac{v^\top R v}{2}+Q(t,x),
\end{equation*}
where $R$ is a positive-definite matrix, $K,\psi,Q$ are smooth and bounded below, and $\rho_0$ is a compactly supported absolutely continuous probability measure with bounded a density. In particular, $z_{a,x}$ is unique for Lebesgue a.e. $x$, and one can choose $z_{a,x}$ in such a way that $(t,x)\mapsto z_{a,x}(t)$ is Borel measurable.

Utilizing the value-function representation \ref{eq:z_ax} of $\phi_a$, we obtain the following saddle-point formluation of \ref{eq:a_optimization}:
\begin{equation}\label{eq:a_primaldual}
\begin{split}
&\inf_a \sup_{z_x:z_x(0)=x} \frac{1}{2}\int_0^T \sum_{i,j=1}^r (\rmK^{-1})_{ij} a_i(t) a_j(t) dt\\
&- \int_{\sR^d}\left[ \int_0^T \left\{L(s,z_x(s),\dot{z}_x(s))+ \sum_{i=1}^r a_i(s)\zeta_i(z_x(s))\right\} ds+\psi(z_x(T)) \right]  d\rho_0(x)
\end{split}
\end{equation}
This saddle-point formulation is the basis of our algorithm in Section \ref{sec:algorithm}.
    
\section{Random Features}\label{sec:rfeatures}

Random features is a simple yet powerful technique to approximate translation invariant positive definite kernels~\cite{RahimiR07}. The foundation of the method is Bochner's theorem from harmonic analysis.

\begin{theorem}[Bochner~\cite{rudin1962fourier}]\label{thm: bochner}
A continuous symmetric shift-invariant kernel $K(x,y) = K(x-y)$ on $\R^d$ is positive definite if and only if $K(\cdot)$ is the Fourier transform of a non-negative measure.
\end{theorem}

Thus, if $K(x-y)$ is a continuous symmetric positive definite kernel there exists a probability distribution $p$ such that
\begin{equation}\label{eq:K_Fourier}
\begin{split}
    K(x-y) =& K(0)~\int_{\sR^d}  e^{i \omega \cdot (x-y)} p(\omega) d\omega\\
    =&K(0)~\mathbb{E}_{\omega \sim p}  \left[ \cos(\omega \cdot x) \cos (\omega \cdot y) + \sin(\omega \cdot x) \sin (\omega \cdot y) \right].
\end{split}
\end{equation}
Hence, we can approximate $K(x-y)$ by sampling $\{\omega_j\}$ iid from $p$:
\begin{equation}\label{eq:K_rfeatures}
    K(x-y)\approx K_r(x-y)=\sum_{j=1}^{r/2} \left[\zeta_j^c(x)\zeta_j^c(y)+\zeta_j^s(x)\zeta_j^s(y)\right],
\end{equation}
where
\begin{equation}\label{eq:zeta_cos_sin}
    \zeta_j^c(x)=\sqrt{\frac{2K(0)}{r}} \cos(\omega_j \cdot x),\quad \zeta_j^s(x)=\sqrt{\frac{2K(0)}{r}} \sin(\omega_j \cdot x).
\end{equation}

Note that this approximation is also shift-invariant, which is a significant advantage for crowd-motion type models where agents interact through their relative positions in the state space. Furthermore, $\rmK$ in \ref{eq:kernel_approx} is the identity matrix, which leads to simple update rules for nonlocal variables $a_1,a_2,\cdots,a_r$: see Section \ref{sec:algorithm}.

The approximation above is viable if one can efficiently sample from $p$. In this paper, we consider Gaussian no-collision repulsive kernels similar to those in~\cite{onken2021multi,onken2021neural}:
\begin{equation}\label{eq:K_Gaussian}
    K(x-y)=\mu \exp \left( -\frac{\|x-y\|^2}{2\sigma^2}\right).
\end{equation}
In this case, one can easily sample from $p$ because it is a Gaussian normal distribution:
\begin{equation}\label{eq:p_Gaussian}
    p(\omega)=\frac{\sigma^d}{(2\pi)^{\frac{d}{2}}} \exp \left(-\frac{\sigma^2 \|\omega\|^2}{2} \right).
\end{equation}

    \section{Trajectory Generation}\label{sec:algorithm}
    
    Here, we propose a primal-dual algorithm inspired by~\cite{nursaude18} to solve \ref{eq:a_primaldual}. Note that the $\sup$ part of \ref{eq:a_primaldual} is a classical optimal control or trajectory optimization problem where the dual variable $a=(a_i)$ acts as a parameter. Thus, we successively optimize trajectories and update the dual variable.
    
    While there exist many trajectory optimization methods for~\ref{eq:a_primaldual}~\cite{nakamura2019adaptive,parkinson2020model,ruthotto2020machine,onken2021neural,onken2021ot,onken2021neural}, we use the direct transcription approach for simplicity~\cite{enright1992discrete}.
    The direct transcription approximates the solution to~\ref{eq:a_primaldual} by discretizing the trajectories over time using, for instance, Euler's Method for the ODE and a midpoint rule to discretize the time integral. 
    Consider a uniform time discretization $0=t_1 < t_2 < \ldots < t_N = T$, and denote the discretized states by $\mathbf{z} = (z(t_1), z(t_2), \ldots, z(t_N))$ and the discretized dual variables by $\mathbf{a}_i = (a_i(t_1), a_i(t_2), \ldots, a_i(t_N))$. The direct transcription approach solves the discretized problem given by 
    \begin{equation}
    \label{eq:a_primaldual_discretized}
    \begin{split}
        &\inf_{\mathbf{a}} \sup_{\mathbf{v}} \frac{1}{2} \sum_{l=1}^{N} h \sum_{i,j=1}^r (\rmK^{-1})_{ij} \mathbf{a}_i[l] \mathbf{a}_j[l]\\
        &- \frac{1}{M}\sum_{m=1}^{M}\left[ \sum_{l=1}^{N} h\left[L(t_l, \mathbf{z}^{\rvv}_{x_m}[l],\mathbf{v}_{x_m}[l])+ \sum_{i=1}^r \mathbf{a}_i[l]\zeta_i(\mathbf{z}^{\rvv}_{x_m}[l])\right] +\psi(\mathbf{z}^{\rvv}_{x_m}[N]) \right], 
    \end{split}
    \end{equation}
    where $\rvv$ is the discretized control, and 
    \begin{equation*}
        \rvz^{\rvv}_{x_m}[l+1]=\mathbf{z}^\rvv_{x_m}[l] + h \mathbf{v}_{x_m}[l], \quad \mathbf{z}^\rvv_{x_m}[1] = x_m,\quad 1\leq l \leq N-1.
    \end{equation*}
    Thus, $\mathbf{z}^{\rvv}_{x_m}[l]$ is the value of $z(t_l)$ for the initial condition $x_m$ and control $v=\dot{z}$. Here, $x_1, x_2, \ldots, x_m$ are samples of initial conditions drawn from $\rho_0$.
    The inner sup problem occurs over the discretized controls 
    \begin{equation}
        \mathbf{v} = 
        \begin{bmatrix}
        \mathbf{v}_{x_1}[1] &  \mathbf{v}_{x_1}[2] & \ldots & \mathbf{v}_{x_1}[N]
        \\
        \mathbf{v}_{x_2}[1] &  \mathbf{v}_{x_2}[2] & \ldots & \mathbf{v}_{x_2}[N]
        \\
        \vdots
        \\
        \mathbf{v}_{x_M}[1] &  \mathbf{v}_{x_M}[2] & \ldots & \mathbf{v}_{x_M}[N]
        \end{bmatrix},
    \end{equation}
    where each row represents the controls for the trajectory defined by initial condition $x_m$.
    The outer inf problem occurs over the discretized coefficients
    \begin{equation}
        \mathbf{a} = 
        \begin{bmatrix}
        \mathbf{a}_1[1] &  \mathbf{a}_1[2] & \ldots & \mathbf{a}_1[N]
        \\
        \mathbf{a}_2[1] &  \mathbf{a}_2[2] & \ldots & \mathbf{a}_2[N]
        \\
        \vdots
        \\
        \mathbf{a}_r[1] &  \mathbf{a}_r[2] & \ldots & \mathbf{a}_r[N]
        \end{bmatrix}.
    \end{equation}
    Indeed, any optimization algorithm can be used to solve this problem. While we used an Euler discretization of the dynamics, any other method could also be used, e.g., RK4. As in~\cite{nursaude18}, we use a version of primal-dual hybrid gradient (PDHG) algorithm~\cite{champock11} to approximate the solution to~\ref{eq:a_primaldual_discretized}. Denoting by
    \begin{equation*}
    \begin{split}
        &\mathcal{L}(\rva, \rvv)= \frac{1}{2} \sum_{l=1}^{N} h \sum_{i,j=1}^r (\rmK^{-1})_{ij} \mathbf{a}_i[l] \mathbf{a}_j[l]\\
        &- \frac{1}{M}\sum_{m=1}^{M}\left[ \sum_{l=1}^{N} h\left[L(s, \mathbf{z}_{x_m}[l],\mathbf{v}_{x_m}[l])+ \sum_{i=1}^r \mathbf{a}_i^{(k)}[l]\zeta_i(\mathbf{z}_{x_m}[l])\right] +\psi(\mathbf{z}_{x_m}[N]) \right],
    \end{split}
    \end{equation*}
    \ref{eq:a_primaldual_discretized} reduces to
    \begin{equation*}
        \inf_{\rva} \sup_{\rvv} \mathcal{L}(\rva,\rvv),
    \end{equation*}
    and the algorithm successively performs the updates
    \begin{equation}\label{eq:algo}
        \begin{split}
            &\mathbf{v}^{(k+1)} = \mathbf{v}^{(k)} + h_v~\nabla_{\rvv} \mathcal{L}(\rva^k, \rvv),
        \\
        &\overline{\rvv}^{(k+1)} = 2\rvv^{(k+1)} - \rvv^{(k)}, 
        \\ 
        &\mathbf{a}^{(k+1)} = \argmin_{\rva} \mathcal{L}(\rva,\overline{\rvv}^{(k+1)})+\frac{\|\rva-\rva^k\|^2}{2h_a},
        \end{split}
    \end{equation}
    where $h_v,h_a>0$ are suitably chosen time-steps, and $(\rva^{(0)},\rvv^{(0)})$ are chosen randomly.
    \begin{remark}
    In the original PDHG of Chambolle and Pock~\cite{champock11} the coupling between $\rva,\rvv$ is bilinear, $\mathcal{L}$ is concave in $\rvv$, and the gradient ascent step in $\rvv$ is replaced by a proximal step. Despite these differences, \ref{eq:algo} has a reliable performance.
    \end{remark}
    
    The gradient ascent step in $\rvv$ is implemented via back-propagation, whereas the proximal step in $\rva$ admits a closed-form solution
    \begin{equation*}
        \begin{split}
        &\mathbf{a}^{(k+1)}[l] = (\mathbf{I} - h_a\mathbf{K}^{-1}) \mathbf{a}^{(k)}[l] + \frac{h_a}{M} 
        \begin{bmatrix}
            \sum_{m=1}^{M} \zeta_1(\mathbf{z}^{\overline{{\rvv}}^{k+1}}_{x_m}[l])
            \\
            \sum_{m=1}^{M} \zeta_2(\mathbf{z}^{\overline{{\rvv}}^{k+1}}_{x_m}[l])
            \\
            \vdots
            \\
            \sum_{m=1}^{M} \zeta_r(\mathbf{z}^{\overline{{\rvv}}^{k+1}}_{x_m}[l])
        \end{bmatrix}, \quad 1 \leq l \leq N.
        \end{split}
    \end{equation*}
    
    Note that for $m \neq m'$ the updates of $\rvv_{x_m}$ and $\rvv_{x_{m'}}$ are decoupled within the $\rvv$ update because the coupling variable $\rva$ is fixed within this update. Furthermore, the random-features approximation yields $\rmK=\operatorname{Id}$, which leads to extremely simple proximal updates for $\rva$:
    \begin{equation*}
        \rva^{(k+1)}_i[l]=(1-h_a)\rva^{(k)}_i[l]+h_a \frac{\sum_{m=1}^{M} \zeta_i(\mathbf{z}^{\overline{{\rvv}}^{k+1}}_{x_m}[l])}{M}, \quad 1 \leq l \leq N.
    \end{equation*}

\section{Numerical Experiments}\label{sec:numerics}

We discuss several numerical examples to demonstrate the efficiency and robustness of our algorithm. The experiments are organized in three groups, A, B, and C, which are presented in Sections \ref{subsec:exA}, \ref{subsec:exB}, and \ref{subsec:exC}, respectively. In experiments A and B we consider high-dimensional problems with low-dimensional interactions - this setting is realistic in the physical setting, e.g., controlling swarm UAVs, since it is often the case that one may have a high-dimensional state/control but the interaction only occurs in the spatial dimensions. In experiment C we consider high-dimensional problems with high-dimensional interactions. The experiments are performed in $d=2,50,100$ dimensions with a fixed time horizon $T=1$. 

\subsection{Experiment A}\label{subsec:exA}

We assume that agents are initially distributed according to a mixture of eight Gaussian distributions centered at the vertices of a regular planar octagon. More precisely, we suppose that
\begin{equation}\label{eq:rho0_A}
    \rho_0(x) \propto \sum^8_{j=1} \exp{ \left( - \frac{\|x-y_i\|^2}{2 \cdot 0.1^2}\right)}
\end{equation}
where
\begin{equation*}
    y_j = \left(\cos\left(\frac{2\pi j}{8}\right), \: \sin\left(\frac{2\pi j}{8}\right), 0, \cdots, 0\right) \in \R^d,\quad 1\leq j \leq 8.
\end{equation*}
Furthermore, we assume that the interaction kernel has the form
\begin{equation}\label{eq:kernel_low_d}
    K(x-y)=\mu \exp \left( -\frac{\|x'-y'\|^2}{2\sigma^2}\right),\quad x,y \in \R^d,
\end{equation}
where $x'=(x_1,x_2) \in \R^2$ for $x=(x_1,x_2,\cdots,x_d) \in \R^d$. Such kernels are repulsive, where $\mu$ is the repulsion intensity, and $\sigma$ is the repulsion radius. Thus, larger $\mu$ leads to more crowd averse agents. Furthermore, the smaller $\sigma$ the more sensitive are the agents to their immediate neighbors. Hence, $\sigma$ can also be interpreted as a safety radius for collision-avoidance applications~\cite{onken2021multi,onken2021neural}. For experiments in A we take $\mu=10$, and $\sigma=0.2,1.25$.

\begin{figure}[h!]
        \centering
         \begin{subfigure}[t]{0.70\textwidth}
             \centering
             \includegraphics[width=\textwidth]{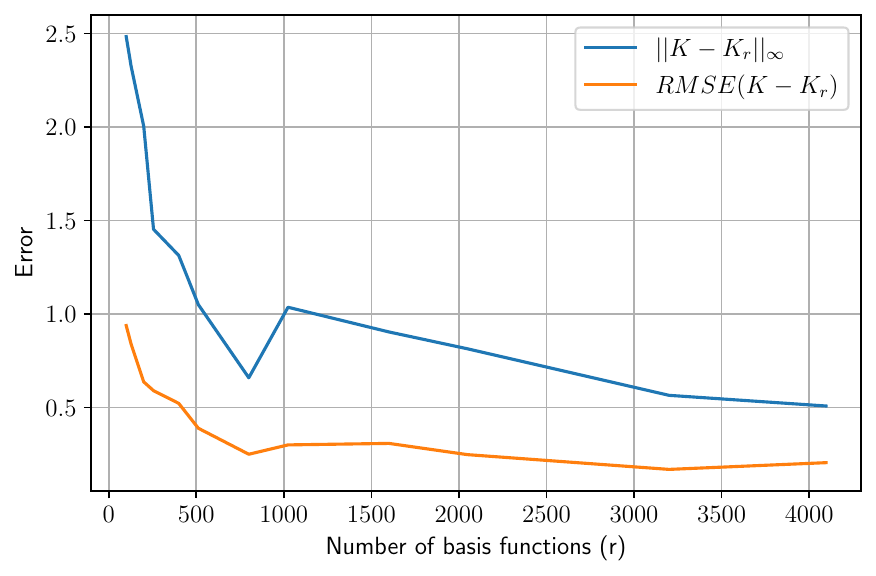}
             \caption{Convergence of errors.}
             \label{fig:kernel-error-convergence-sigma-02}
         \end{subfigure} \\
         
         \begin{subfigure}[t]{0.49\textwidth}
             \centering
             \includegraphics[width=\textwidth]{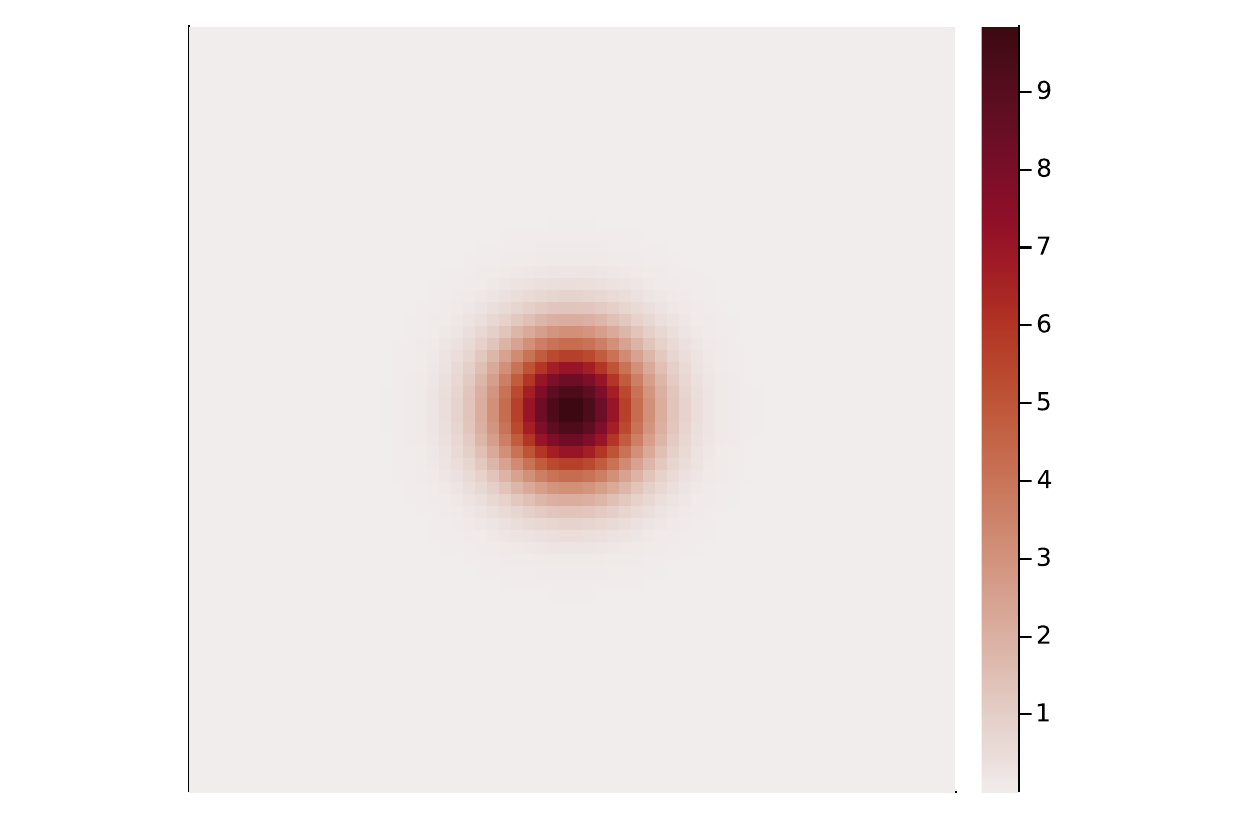}
             \caption{True Kernel}
             \label{fig:true-kernel-2d-plot-sigma-02}
         \end{subfigure}
         \begin{subfigure}[t]{0.49\textwidth}
             \centering
             \includegraphics[width=\textwidth]{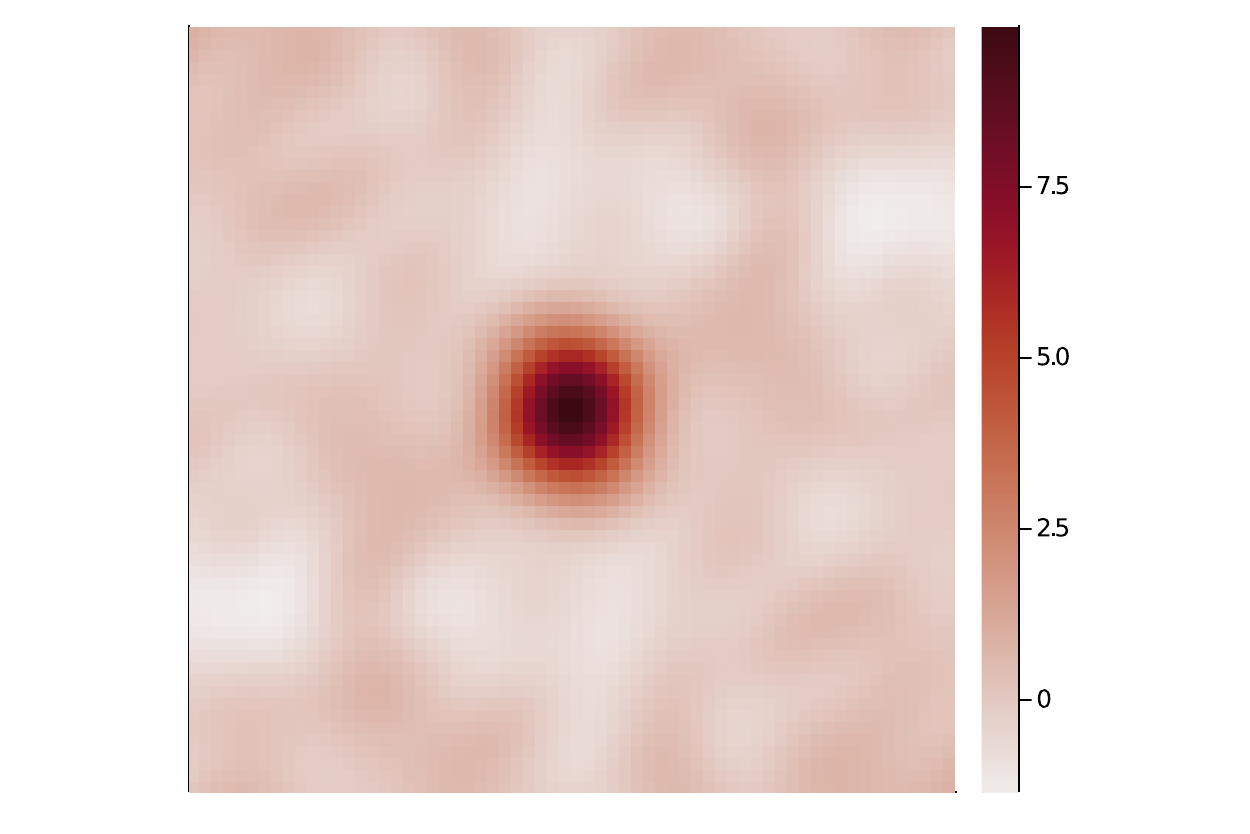}
             \caption{Approximate kernel with \\ $r=512$ features.}
             \label{fig:approx-kernel-2d-plot-sigma-02}
         \end{subfigure}
        \caption{Kernel approximation for $\sigma = 0.2, \mu = 10.0$.}
        \label{fig:kernel-approx-sigma-02}
    \end{figure}

\begin{figure}[h!]
        \centering
         \begin{subfigure}[t]{0.70\textwidth}
             \centering
             \includegraphics[width=\textwidth]{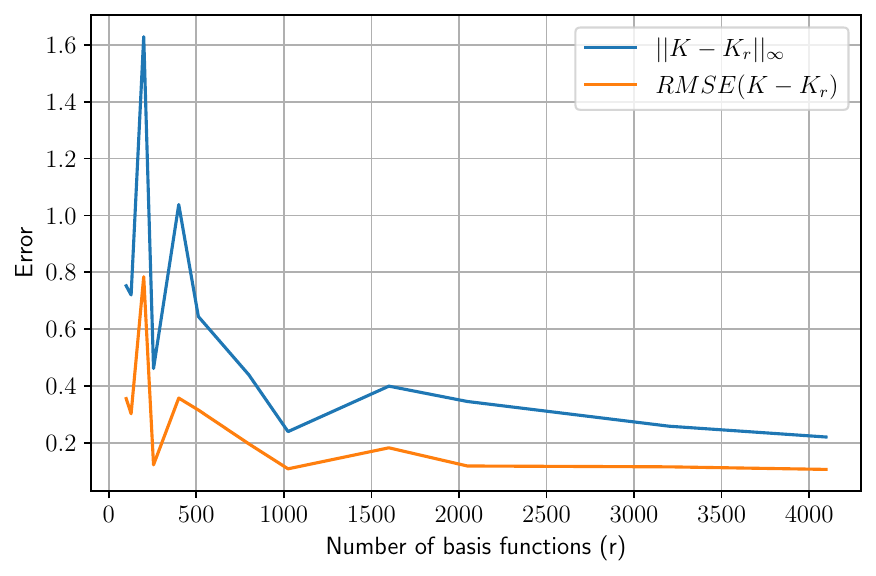}
             \caption{Convergence of errors}
             \label{fig:kernel-error-convergence-sigma-125}
         \end{subfigure} \\
         
         \begin{subfigure}[t]{0.49\textwidth}
             \centering
             \includegraphics[width=\textwidth]{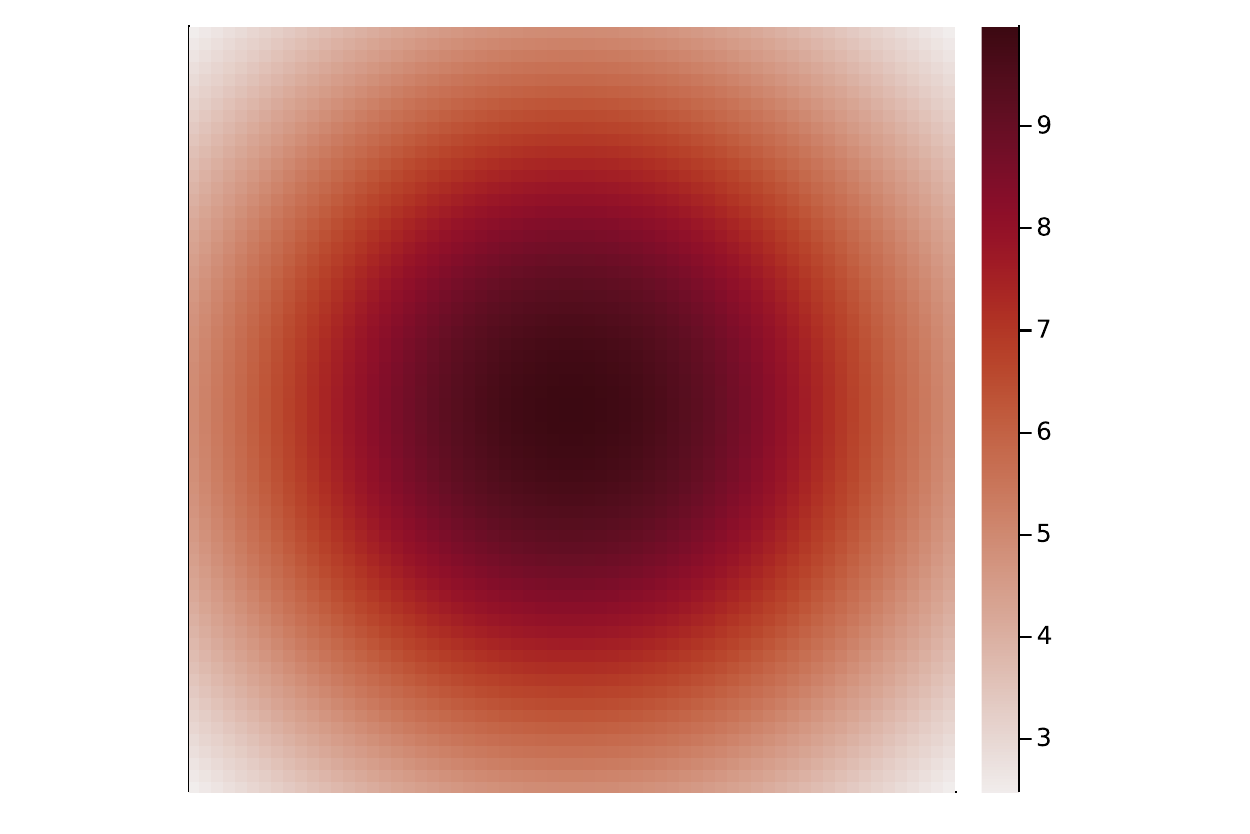}
             \caption{True Kernel}
             \label{fig:true-kernel-2d-plot-sigma-125}
         \end{subfigure}
        \begin{subfigure}[t]{0.49\textwidth}
             \centering
             \includegraphics[width=\textwidth]{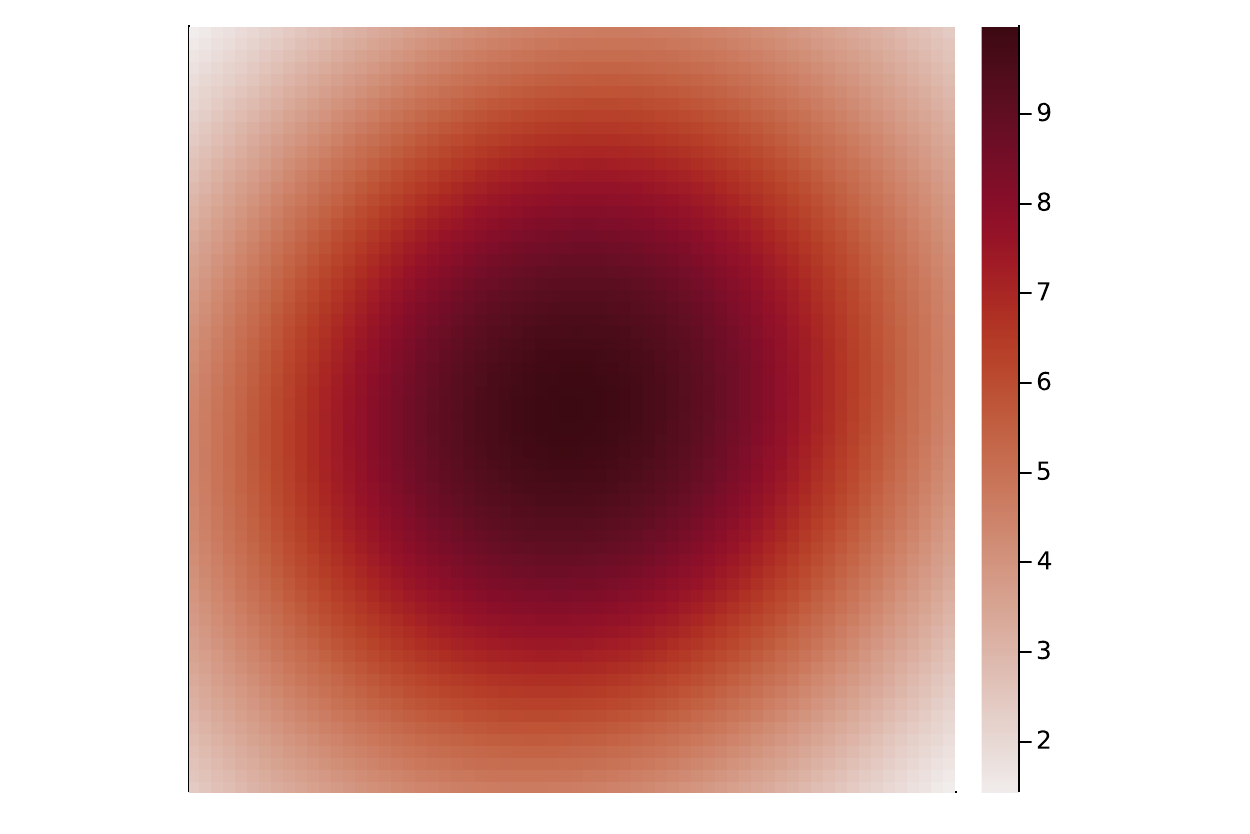}
             \caption{Approximate kernel with  \\ $r=512$ features.}
             \label{fig:approx-kernel-2d-plot-sigma-125}
         \end{subfigure}
        \caption{Kernel approximation for $\sigma = 1.25, \mu = 10.0$.}
        \label{fig:kernel-approx-sigma-125}
    \end{figure}

The random-features approximation of $K$ is given by
\begin{equation*}
K(x-y)\approx K_r(x-y)=\sum_{j=1}^{r/2}\left[\zeta_j^c(x')\zeta_j^c(y')+\zeta_j^s(x') \zeta_j^s(y') \right],
\end{equation*}
where
\begin{equation*}
    \zeta_j^c(x')=\sqrt{\frac{2\mu}{r}} \cos(\omega'_j \cdot x'),\quad \zeta_j^s(x')=\sqrt{\frac{2\mu}{r}} \sin(\omega'_j \cdot x'),\quad x'\in \R^2,
\end{equation*}
and $\{\omega'_j\}_{j=1}^{r/2} \subset \R^2$ are drawn randomly from
\begin{equation*}
    p(\omega')=\frac{\sigma^2}{2\pi} \exp \left(-\frac{\sigma^2 \|\omega'\|^2}{2} \right),\quad \omega' \in \R^2.
\end{equation*}
We plot the convergence of approximate kernels to the true one in Figures \ref{fig:kernel-error-convergence-sigma-02} and \ref{fig:kernel-error-convergence-sigma-125} for $\sigma=0.2$ and $\sigma=1.25$, respectively. This is done by comparing the values generated by the true and approximate kernels $K(x',0)=K(x'-0)$, $K_r(x',0)=K_r(x'-0)$ in $l^\infty$ and $l^2$ norms for $x'$ on a 2-dimensional grid centred at the origin. Further, in Figures \ref{fig:true-kernel-2d-plot-sigma-02}, \ref{fig:approx-kernel-2d-plot-sigma-02} and Figures \ref{fig:true-kernel-2d-plot-sigma-125}, \ref{fig:approx-kernel-2d-plot-sigma-125}, we visually compare the approximation to the true kernel on this grid. 
In experiments A, we choose $r=512$ for both values of $\sigma$. 


We take the Lagrangian and terminal cost functions
\begin{equation*}
    L(t,x,v)=\frac{\|v\|^2}{2},\quad \psi(x)=10 \|x-x_{\text{target}}\|^2,\quad (t,x,v) \in (0,1)\times \R^d \times \R^d,
\end{equation*}
where $x_{\text{target}}=0$. This choice corresponds to a model where crowd-averse agents travel from initial positions towards a target location, $x_{\text{target}}$. Finally, we sample $M=256$ initial positions from $\rho_0$. 

\begin{figure}[h!]
        \centering
         \begin{subfigure}[b]{0.32\textwidth}
             \centering
             \includegraphics[width=\textwidth]{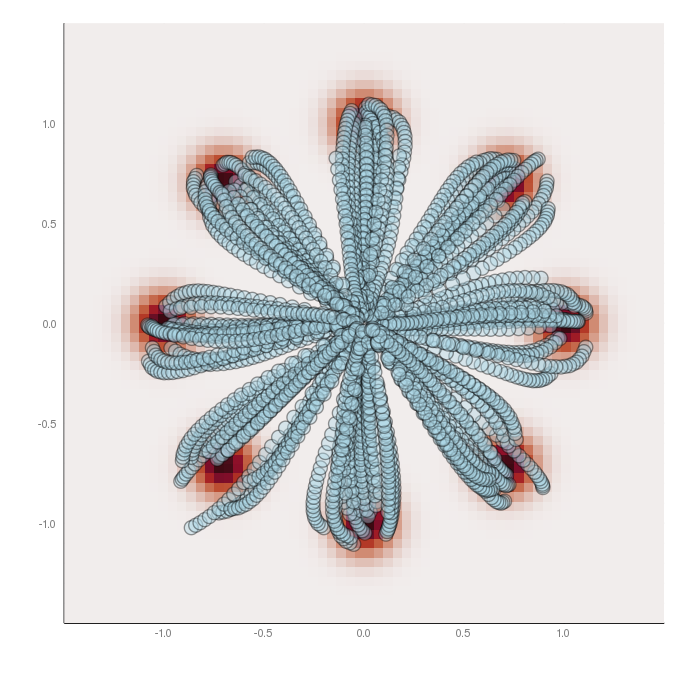}
             \caption{$d=2, \sigma=0.2$.}
             \label{fig:8gaussian-d-2-sigma-02}
         \end{subfigure}
         \begin{subfigure}[b]{0.32\textwidth}
             \centering
             \includegraphics[width=\textwidth]{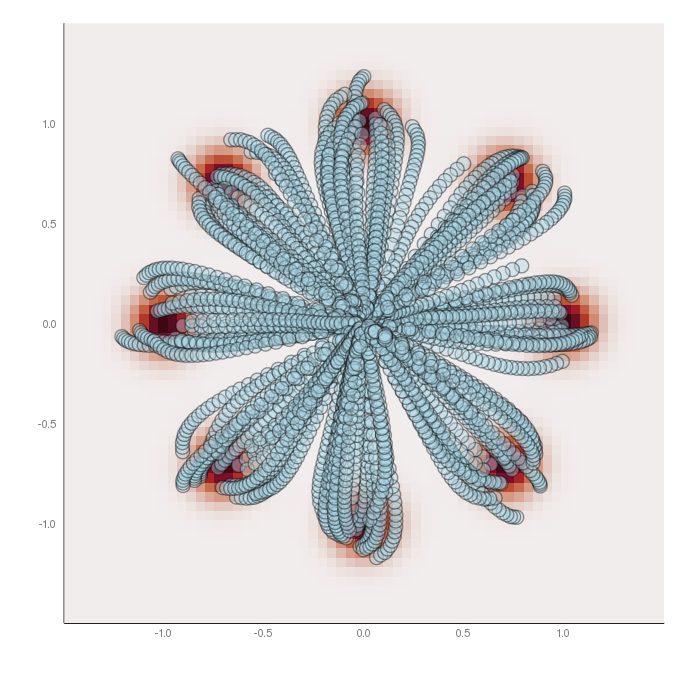}
             \caption{$d=50, \sigma=0.2$.}
             \label{fig:8gaussian-d-50-sigma-02}
         \end{subfigure}
         \begin{subfigure}[b]{0.32\textwidth}
             \centering
             \includegraphics[width=\textwidth]{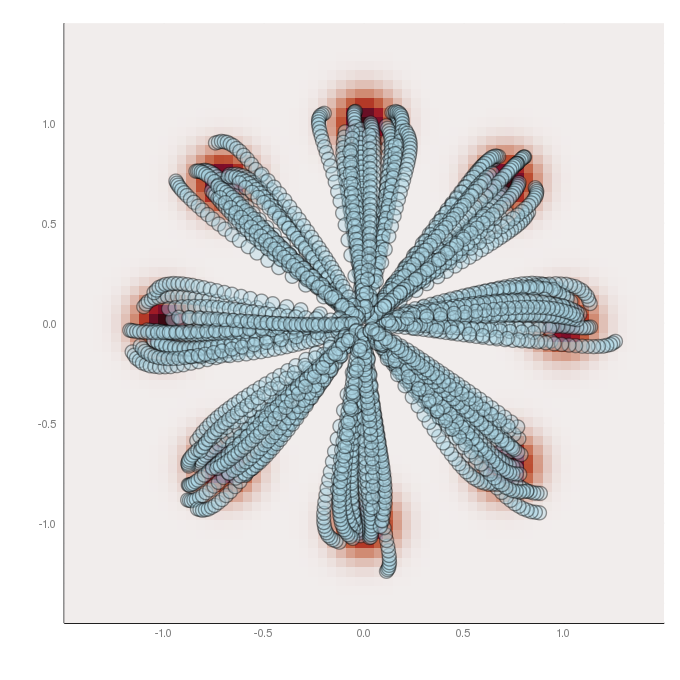}
             \caption{$d=100, \sigma=0.2$.}
             \label{fig:8gaussian-d-100-sigma-02}
         \end{subfigure}
        \caption{Agents' trajectories in experiments A for $\sigma=0.2$ plotted on the first two dimensions. Agents move from 8 Gaussian distributions (colored red) to the target point ($0,0$). Each plot shows the trajectories solved in different dimensions: (a) $d=2$, (b) $d=50$, (c) $d=100$.}
        \label{fig:8gaussian-first-two-02}
    \end{figure}

\begin{figure}[h!]
        \centering
         \begin{subfigure}[b]{0.32\textwidth}
             \centering
             \includegraphics[width=\textwidth]{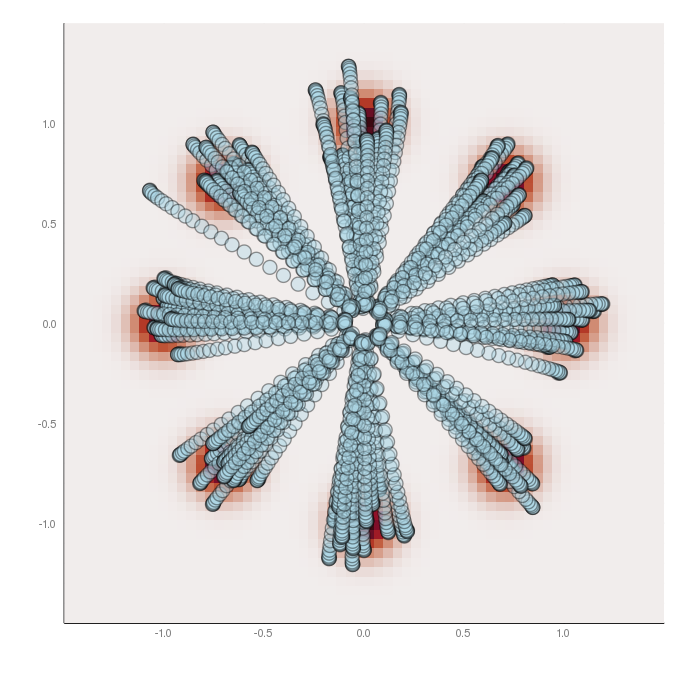}
             \caption*{$d=2, \sigma=1.25$.}
             \label{fig:8gaussian-d-2-sigma-125}
         \end{subfigure}
         \begin{subfigure}[b]{0.32\textwidth}
             \centering
             \includegraphics[width=\textwidth]{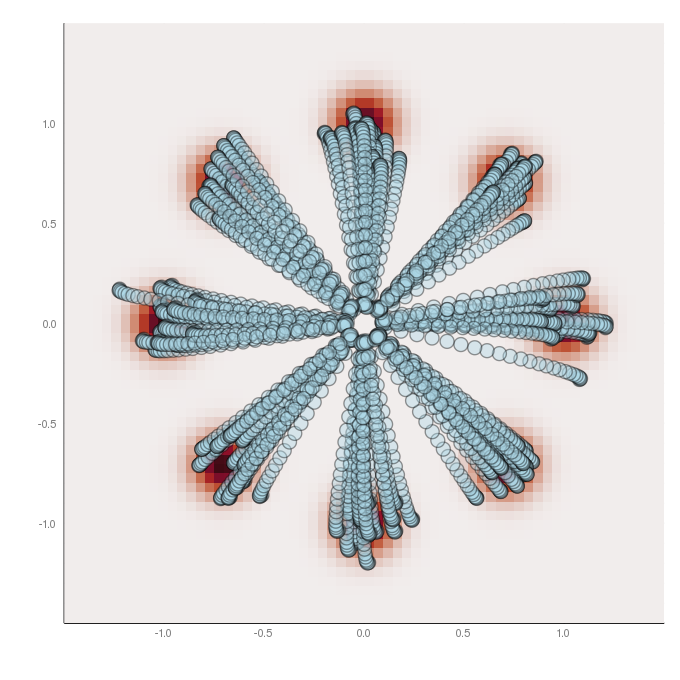}
             \caption{$d=50, \sigma=1.25$.}
             \label{fig:8gaussian-d-50-sigma-125}
         \end{subfigure}
         \begin{subfigure}[b]{0.32\textwidth}
             \centering
             \includegraphics[width=\textwidth]{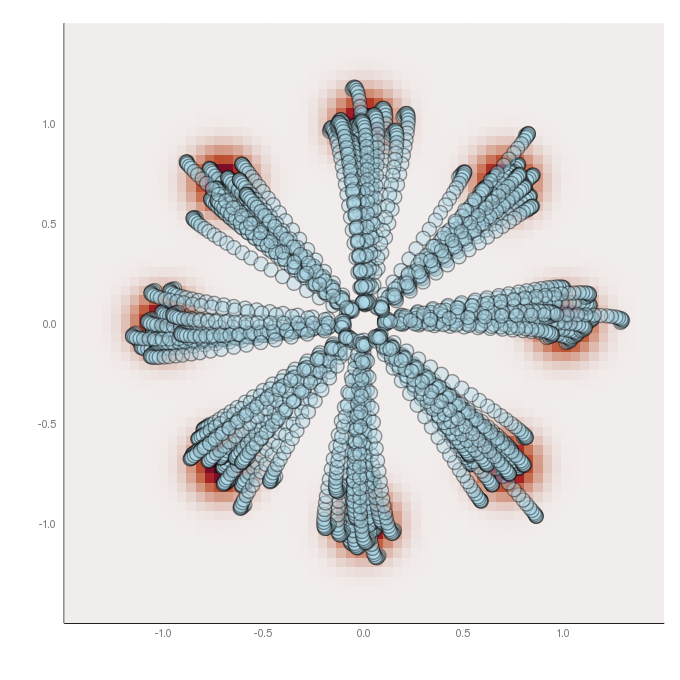}
             \caption{$d=100, \sigma=1.25$.}
             \label{fig:8gaussian-d-100-sigma-125}
         \end{subfigure}
        \caption{Agents' trajectories in experiments A for $\sigma=1.25$ plotted on the first two dimensions. Agents move from 8 Gaussian distributions (colored red) to the target point ($0,0$). Each plot shows the trajectories solved in different dimensions: (a) $d=2$, (b) $d=50$, (c) $d=100$.}
        \label{fig:8gaussian-first-two-125}
    \end{figure}

\begin{figure}[h!]
    \centering
     \begin{subfigure}[b]{0.48\textwidth}
         \centering
         \includegraphics[width=\textwidth]{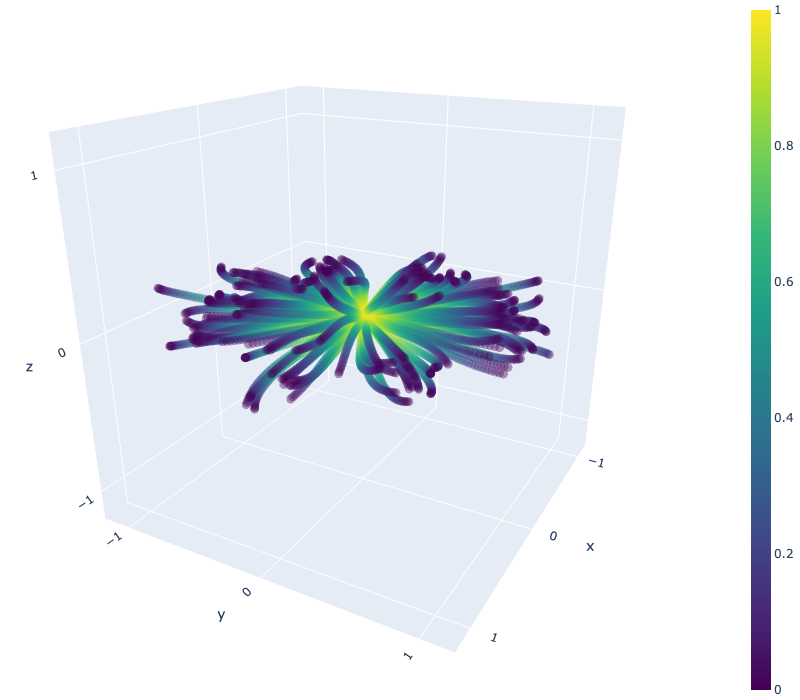}
         \caption{3D view.}
     \end{subfigure}
     \begin{subfigure}[b]{0.48\textwidth}
         \centering
         \includegraphics[width=\textwidth]{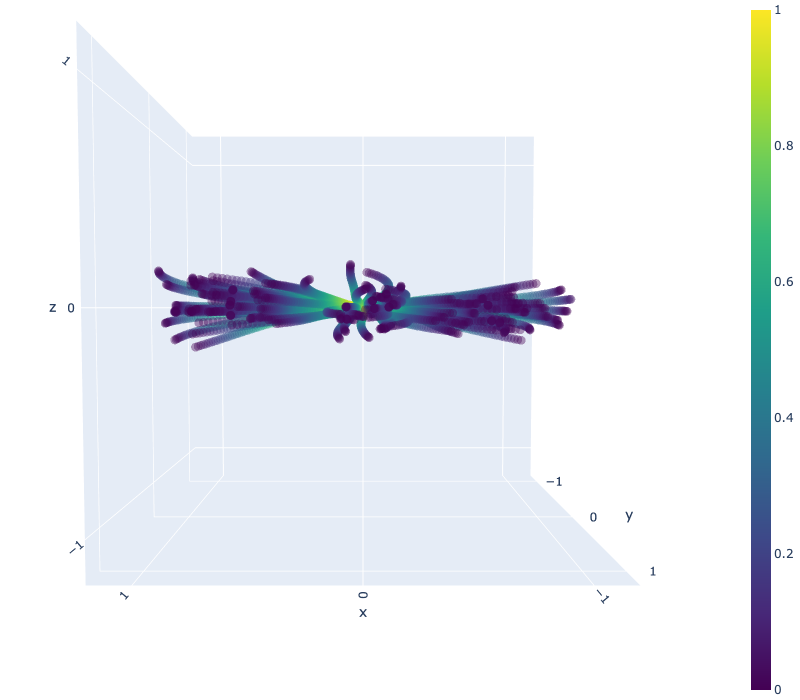}
         \caption{XZ view.}
     \end{subfigure}
     \begin{subfigure}[b]{0.48\textwidth}
         \centering
         \includegraphics[width=\textwidth]{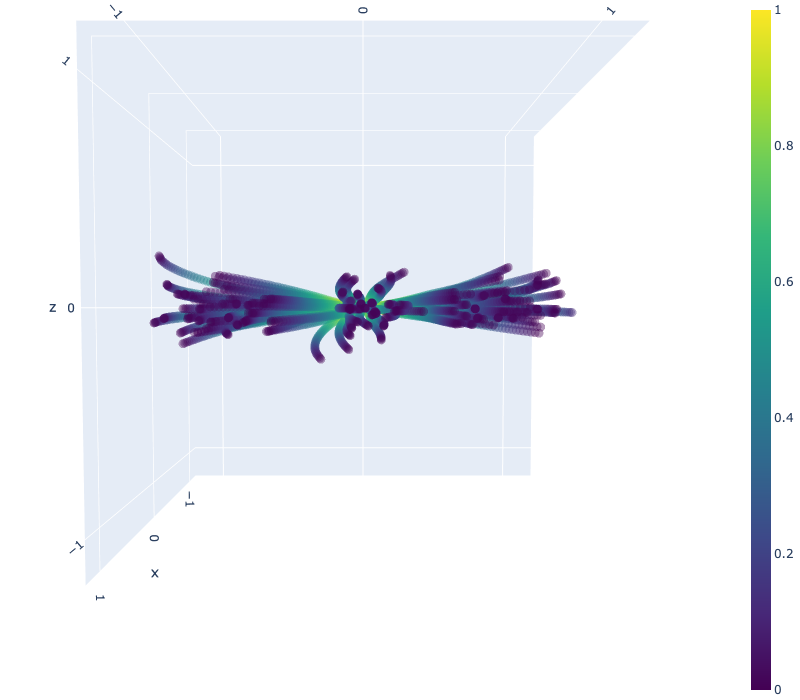}
         \caption{YZ view.}
     \end{subfigure}
     \begin{subfigure}[b]{0.48\textwidth}
         \centering
         \includegraphics[width=\textwidth]{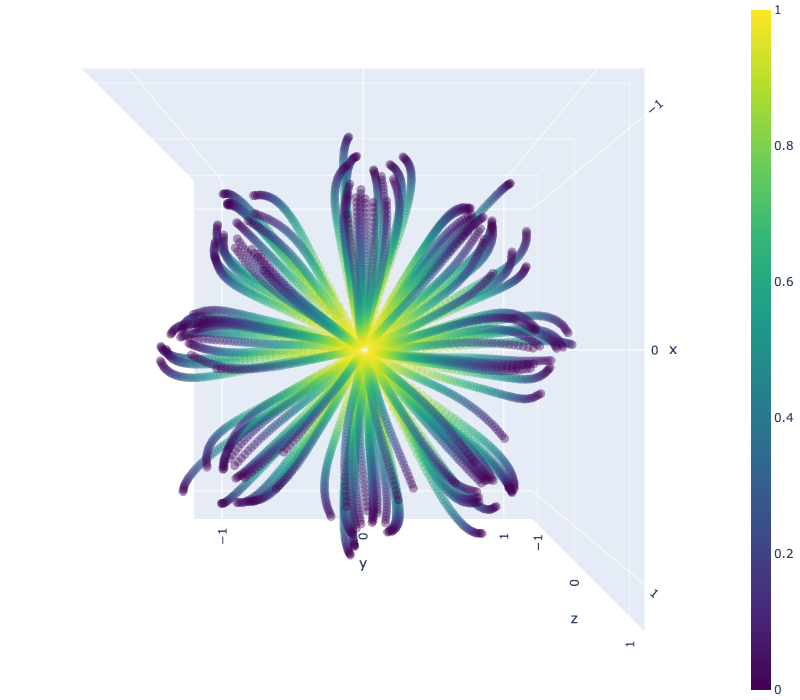}
         \caption{Top view.}
     \end{subfigure}
    \caption{3D plots of agents' trajectories in experiments A with low-dimensional interactions (the first two dimensions) for $d=50$, ${\sigma}=0.2$. The plots show the first three dimensions of the trajectories. Each agent starts from $t=0$ (colored blue) to $t=1$ (colored yellow). The plots are from four different viewpoints. (a): 3D view, (b), (c): side views (XZ view and YZ view), (d): top view (XY view). The interactions of agents are only across the first two dimensions. Thus, while the agents spread in XY axis (see (d) for the top view), they move to the target point almost linearly in other axis (see (b) and (c) for the side views). }
    \label{fig:8gaussian-first-two-3d}
\end{figure}

In Figure \ref{fig:8gaussian-first-two-02} we plot the projections of agents' trajectories on the first two dimensions when the repulsion radius is $\sigma=0.2$ and $d=2,50,100$. Analogously, we plot the agents trajectories for $\sigma=1.25$ and $d=2,50,100$ in Figure \ref{fig:8gaussian-first-two-125}. Note that trajectories split more when $\sigma=0.2$, which corresponds to the case when agents are more sensitive to their immediate neighbors. Additionally, note that the terminal cost function enforces agents to reach the destination $x_{\text{target}}=0$. The 3D trajectories are plotted in Figure \ref{fig:8gaussian-first-two-3d}.

In Table \ref{table:A} we report the population running cost
\begin{equation*}
    \frac{h}{M}\sum_{m=1}^{M} \sum_{l=1}^{N} L(s, \rvz_{x_m}[l],\rvv_{x_m}[l]),
\end{equation*}
interaction cost
\begin{equation*}
    h\sum_{l=1}^{N} \frac{1}{2M^2}\sum_{m,m'=1}^{M} K_r(\rvz_{x_m}[l],\rvz_{x_{m'}}[l])=\frac{h}{2M^2}\sum_{l=1}^{N}\left( \sum_{i=1}^r \zeta_i(\rvz_{x_m}[l])\right)^2,
\end{equation*}
terminal cost
\begin{equation*}
    \frac{1}{M}\sum_{m=1}^{M} \psi(\mathbf{z}_{x_m}[N]),
\end{equation*}
and the total cost at the equilibrium.

\begin{table}[h!]
    \centering
    \begin{tabular}{|c|c|c|c|c|c|}

            \hline
            $d$   & $\sigma$ & Running & Interaction & Terminal & Total\\
            \hline
            2 & 0.2  & 0.526 & 0.465 & 0.0108  & 1.10
            \\
            2 & 1.25  & 0.621 & 3.57 & 0.00997  & 4.29
            \\
            \hline
            50 & 0.2 & 0.754    & 0.454 & 0.0116 &  1.32
            \\
            50 & 1.25 & 0.825    & 3.58 & 0.0109 &  4.51
            \\
            \hline
            100 & 0.2  & 0.992 & 0.533 & 0.0140 & 1.67
            \\
            100 & 1.25  & 1.11 & 3.26 & 0.0139 & 4.51
            \\
            \hline
        \end{tabular}
        \caption{Running, interaction, terminal, and total costs in experiments A.}
    \label{table:A}
\end{table}

\subsection{Experiment B}\label{subsec:exB}

In this set of experiments we assume that agents are initially distributed according to
\begin{equation*}
    \rho_0(x)\propto \exp\left( -\frac{\|x-x_{\text{initial}}\|^2}{2\cdot 0.2^2}\right),\quad x\in \R^d,
\end{equation*}
where $x_{\text{initial}}=(0,1,0,\cdots,0) \in \R^d$. Furthermore, we assume that the Lagrangian and terminal cost functions are
\begin{equation*}
\begin{split}
    L(t,x,v)=&\frac{\|v\|^2}{4}+5\max \left( x'^\top  \begin{bmatrix} 1 & 0 \\ 0 & -5 \end{bmatrix} 
    x', 0
    \right),\quad 
    \psi(x)=10 \|x-x_{\text{target}}\|^2,
\end{split}
\end{equation*}
for $(t,x,v) \in (0,1)\times \R^d \times \R^d$, where $x_{\text{target}}=(0,-1,0,\cdots,0)$.

\begin{figure}[h!]
        \centering
         \begin{subfigure}[t]{0.70\textwidth}
             \centering
             \includegraphics[width=\textwidth]{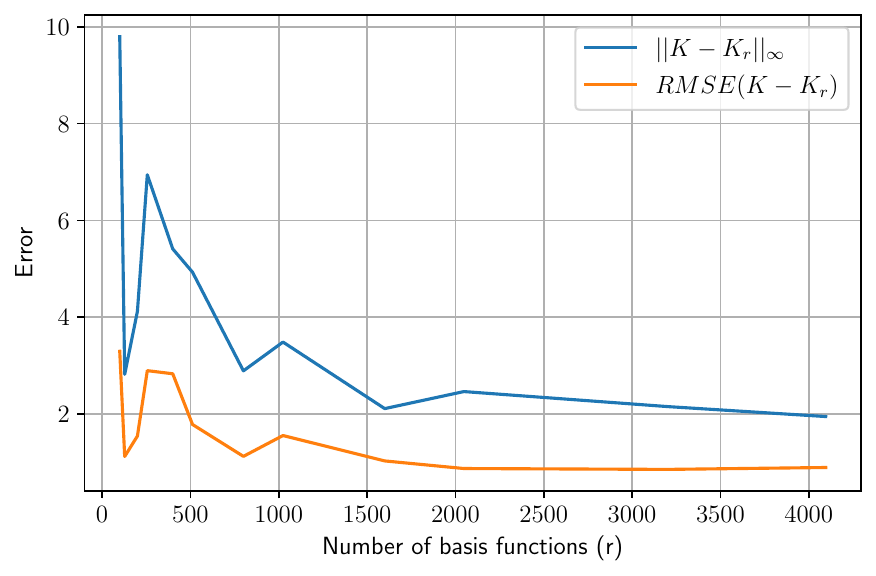}
             \caption{Convergence of errors.}
             \label{fig:kernel-error-convergence-sigma-10}
         \end{subfigure}
         \begin{subfigure}[t]{0.49\textwidth}
             \centering
             \includegraphics[width=\textwidth]{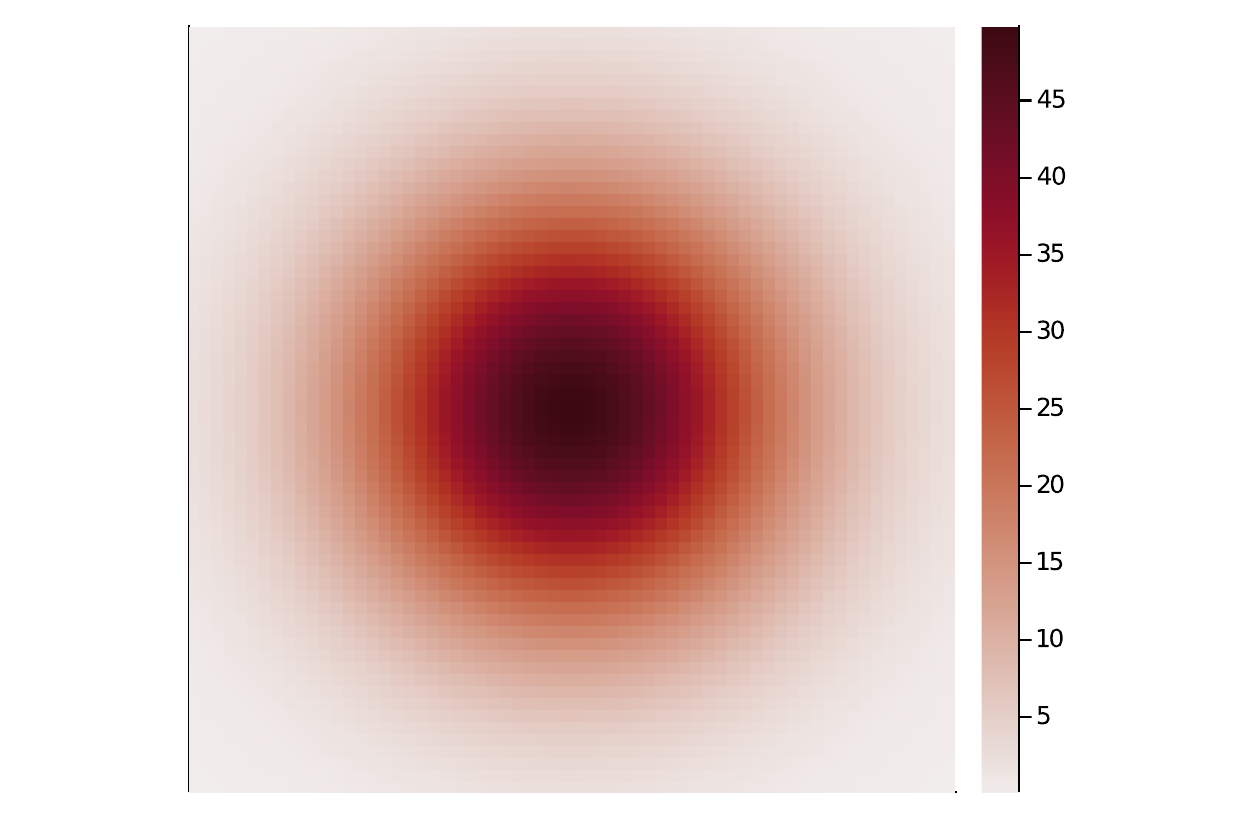}
             \caption{True Kernel}
             \label{fig:true-kernel-2d-plot-sigma-10}
         \end{subfigure}
        \begin{subfigure}[t]{0.49\textwidth}
             \centering
             \includegraphics[width=\textwidth]{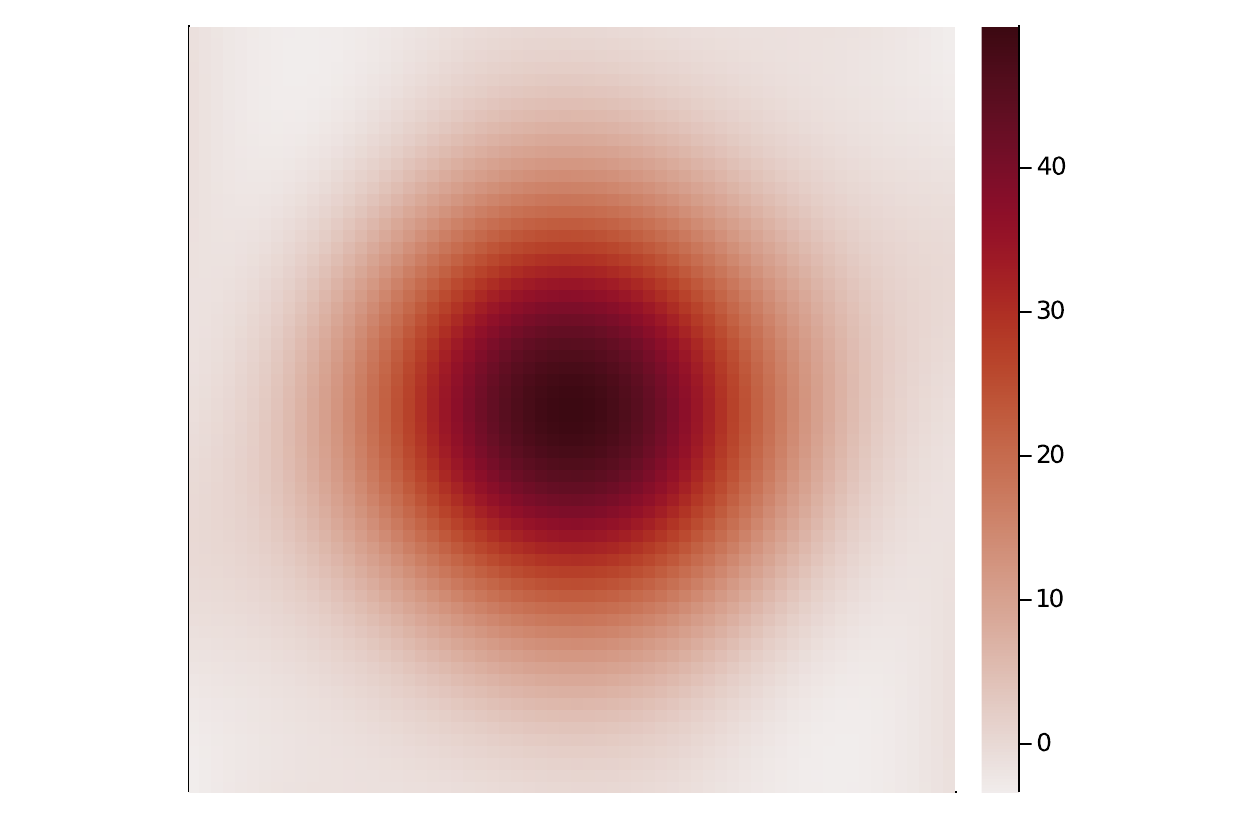}
             \caption{Approximate kernel with \\ $r=512$ features.}
             \label{fig:approx-kernel-2d-plot-sigma-10}
         \end{subfigure}
        \caption{Kernel approximation for $\sigma = 1.0, \mu = 50.0$.}
        \label{fig:kernel-approx-sigma-1}
    \end{figure}

As before, we consider low-dimensional interactions with a kernel of the form \ref{eq:kernel_low_d}. We take $\mu=50$ and $\sigma=1$. The approximation error and the approximate kernel for $r=512$ are plotted in Figure \ref{fig:kernel-approx-sigma-1}. As before, the plots are generated by evaluating the true kernel and the approximate kernel at points on a 2-dimensional grid. 

\begin{figure}[h!]
        \centering
         \begin{subfigure}[b]{0.32\textwidth}
             \centering
             \includegraphics[width=\textwidth]{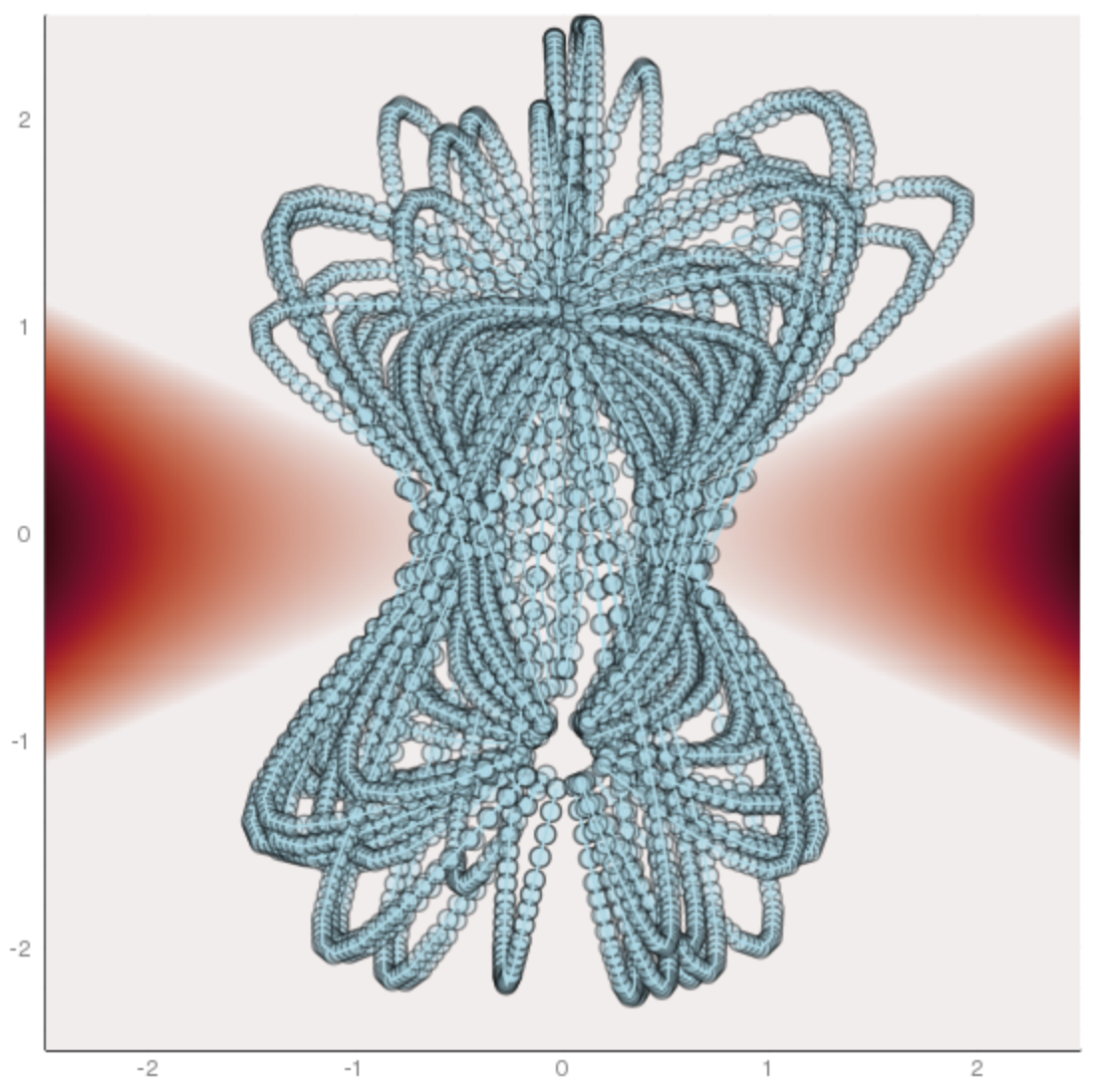}
             \caption{$d=2$.}
             \label{fig:obstacles-d-2-sigma-1}
         \end{subfigure}
         \begin{subfigure}[b]{0.32\textwidth}
             \centering
             \includegraphics[width=\textwidth]{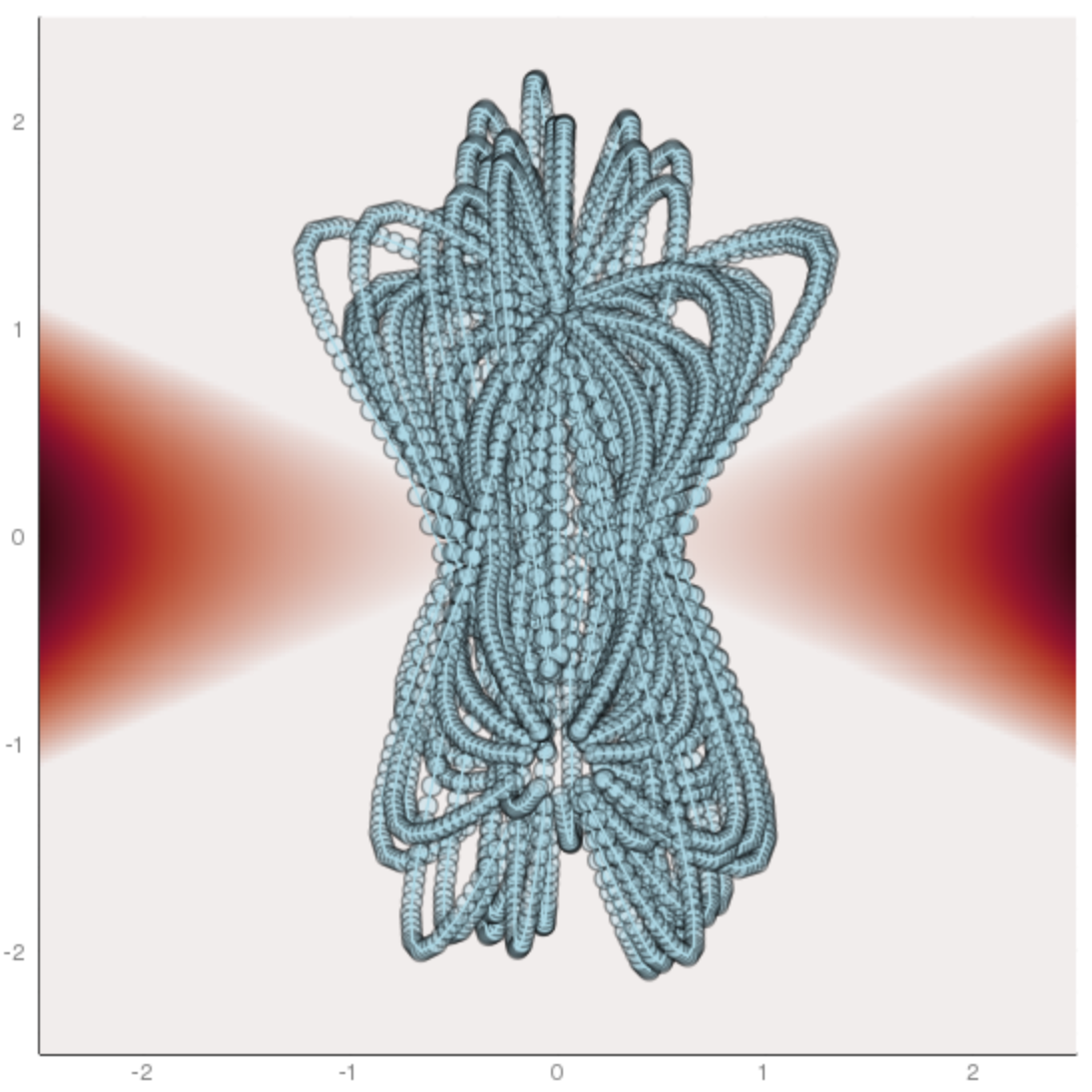}
             \caption{$d=50$.}
             \label{fig:obstacles-d-50-sigma-1}
         \end{subfigure}
         \begin{subfigure}[b]{0.32\textwidth}
             \centering
             \includegraphics[width=\textwidth]{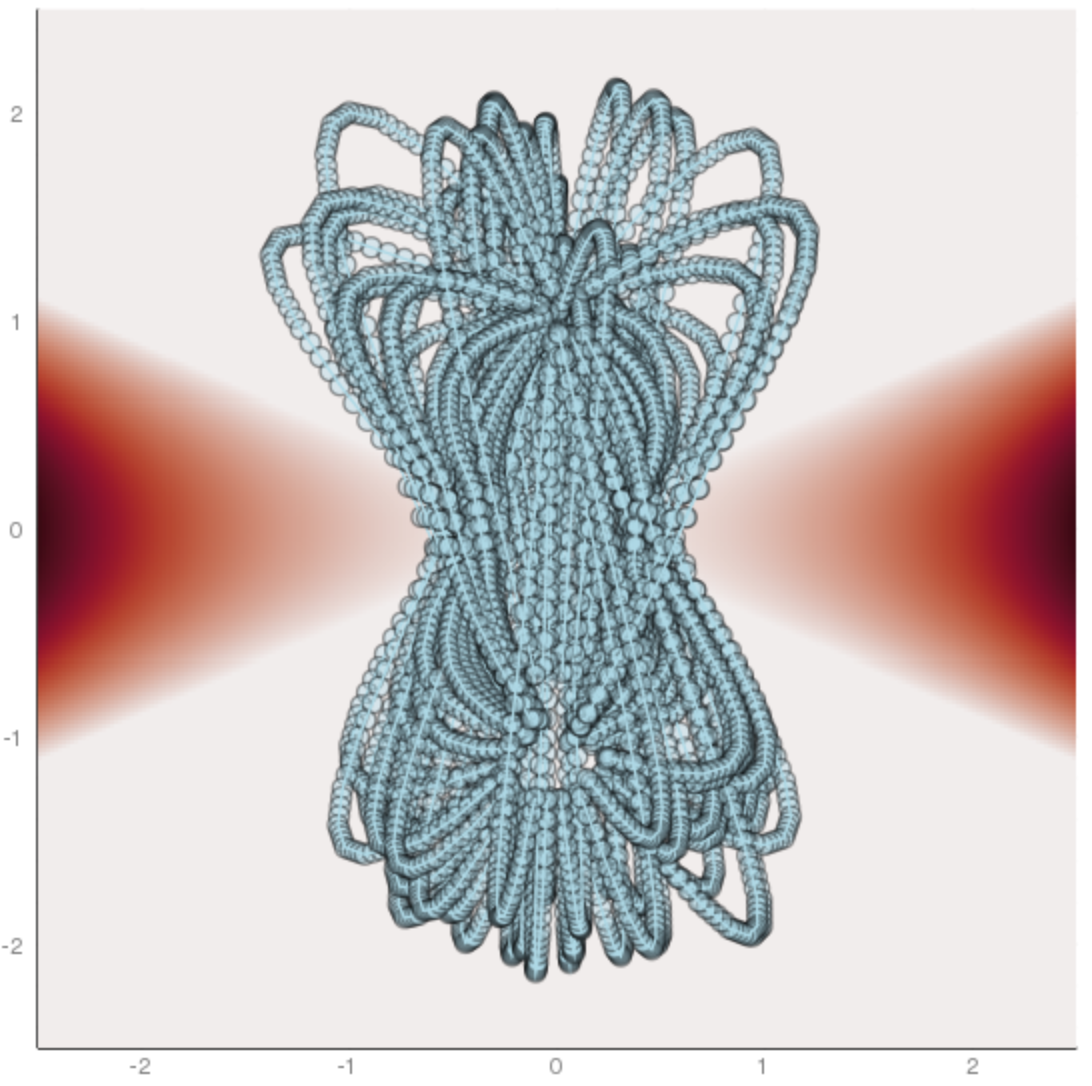}
             \caption{$d=100$.}
             \label{fig:obstacles-d-100-sigma-1}
         \end{subfigure}
        \caption{Agents' trajectories in experiments B plotted on the first two dimensions. Agents move from the initial distribution (near $(0,1)$) to the target point ($0,-1$) while avoiding the obstacle (colored red). Each plot shows the trajectories solved in different dimensions: (a) $d=2$, (b) $d=50$, (c) $d=100$.}
        \label{fig:obstacles}
    \end{figure}

Thus, in experiments B we model a crowd-averse population that travels from around an initial point, $x_{\text{initial}}$, to a target point, $x_{\text{target}}$, avoiding wedge shaped obstacles. The projections of agents' trajectories on the first two dimensions are plotted in Figure \ref{fig:obstacles}.

Note that the trajectories split at close to the initial and target points, demonstrating the crowd-averse behavior of the agents. On the other hand, obstacles force the agents to converge at the bottleneck.

\begin{figure}[h!]
    \centering
     \begin{subfigure}[b]{0.48\textwidth}
         \centering
         \includegraphics[width=\textwidth]{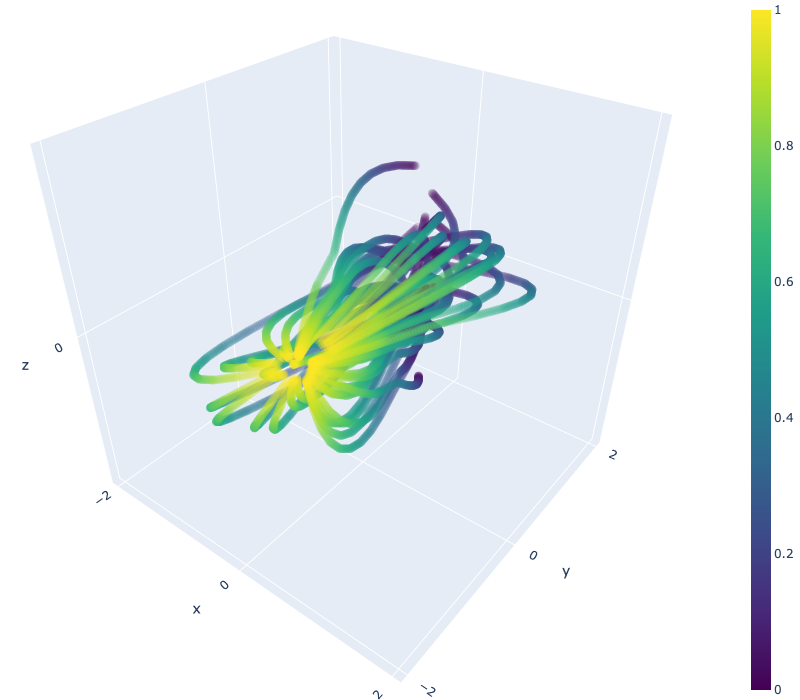}
         \caption{3D view.}
     \end{subfigure}
     \begin{subfigure}[b]{0.48\textwidth}
         \centering
         \includegraphics[width=\textwidth]{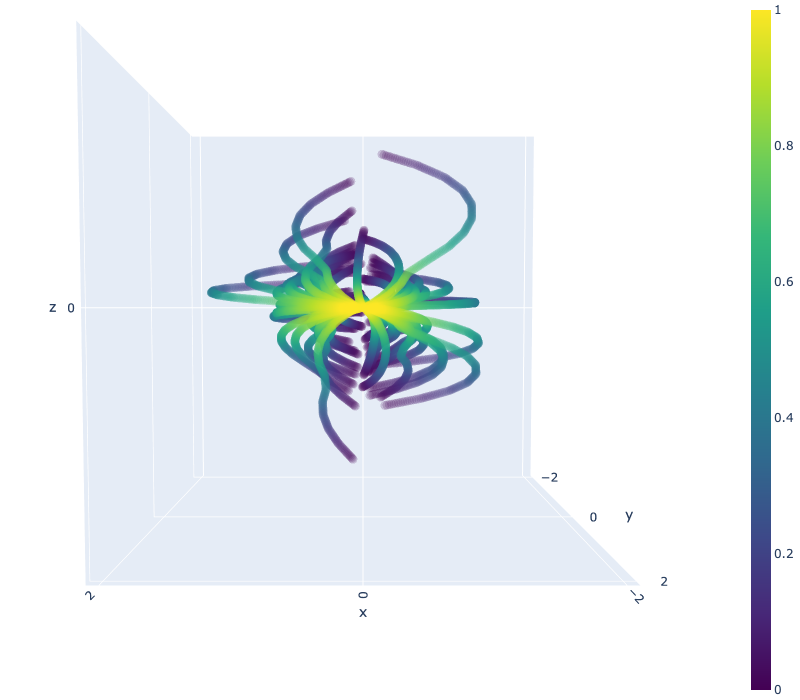}
         \caption{XZ view.}
     \end{subfigure}
     \begin{subfigure}[b]{0.48\textwidth}
         \centering
         \includegraphics[width=\textwidth]{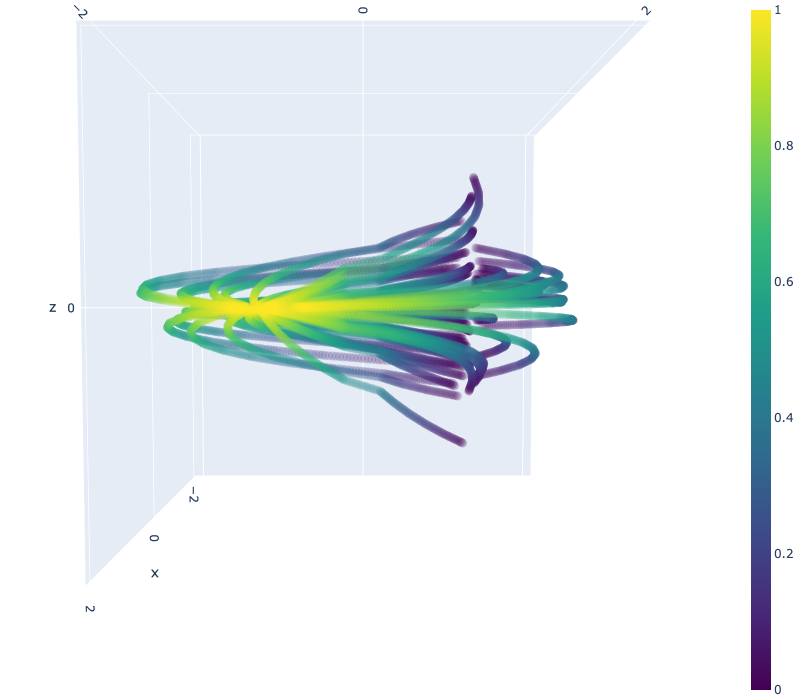}
         \caption{YZ view.}
     \end{subfigure}
     \begin{subfigure}[b]{0.48\textwidth}
         \centering
         \includegraphics[width=\textwidth]{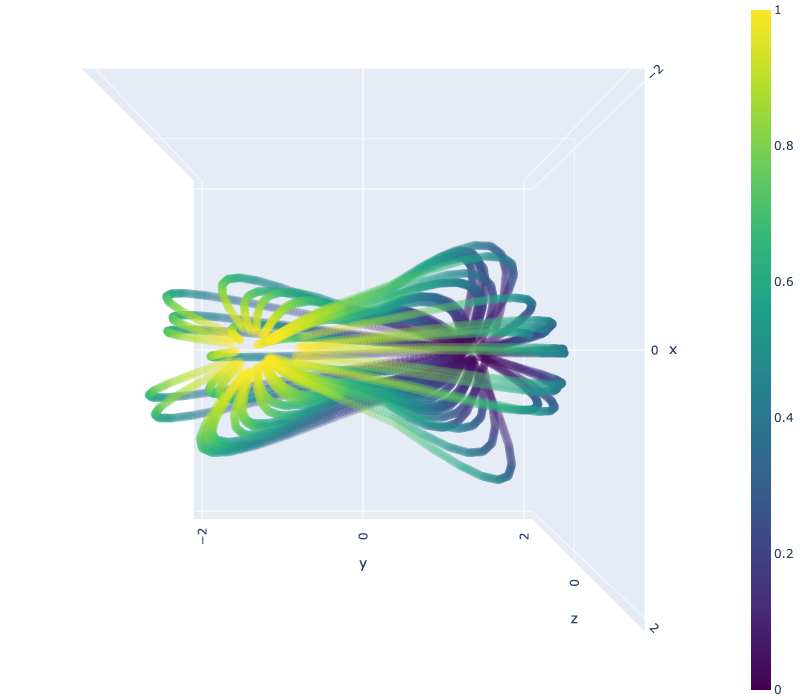}
         \caption{Top view.}
     \end{subfigure}
    \caption{3D plots of Agents' trajectories with low-dimensional interactions in experiments B for $d=50$. The plots show the first three dimensions of the trajectories. Each agent starts from $t=0$ (colored blue) to $t=1$ (colored yellow) while avoiding the obstacle (see Figure~\ref{fig:obstacles}). The plots are from four different viewpoints. (a): 3D view (the target point is at the lower-left side of the plot and the initial distributions are at the top-right side of the plot), (b), (c): side views (XZ view and YZ view), (d): top view (XY view). The interactions of agents are only across the first two dimensions. Thus, while the agents spread in XY axis (see (d) for the top view), they move to the target point almost linearly in other axis (see (c) for the side view). }
    \label{fig:obstacles-all-d-3d}
\end{figure}

We plot the 3D trajectories in Figure \ref{fig:obstacles-all-d-3d} and report running, interaction, terminal, and total costs in Table \ref{table:obstacles}.

\begin{table}[h!]
    \centering
    \begin{tabular}{|c|c|c|c|c|c|}
            \hline
            $d$  & Running & Interaction & Terminal & Total
            \\
            \hline
            2   & 3.72 & 12.2 & 0.388 & 16.3
            \\
            50   & 2.63 & 15.4 & 0.533 & 18.6
            \\
            100   & 2.86 & 14.8 & 0.567 & 18.3
            \\
            \hline
        \end{tabular}
        \caption{Running, interaction, terminal, and total costs in experiments B.}
    \label{table:obstacles}
\end{table}

\subsection{Experiment C}\label{subsec:exC}

In experiments A, B we consider high-dimensional problems with low-dimensional interactions. Here, we perform experiments similar to A but with full-dimensional interactions to demonstrate the efficiency of our method for higher-dimensional interactions as well.

Thus, we assume that we are in the same setup as in A with the only difference that $K$ is a full-dimensional interaction \ref{eq:K_Gaussian} with $\sigma=\hat{\sigma}\cdot \sqrt{d/2}$, and $\mu=10$ and $\mu=1$ for $\hat \sigma=0.2$ and $\hat{\sigma}=1.25$, respectively. Here, $\hat{\sigma}$ is a \textit{dimensionless} repulsion radius. Indeed, since for $\rho_0$ in \ref{eq:rho0_A} the variance of constituent Gaussians is the same across dimensions, the average distance between agents scales with a factor $\sqrt{d}$ near the centers of these Gaussians. Hence, if we used the same repulsion radius across all dimensions, the effective interaction would be different, and it would be hard to interpret the results. By fixing a repulsion radius $\hat{\sigma}$ for $d=2$ and scaling it accordingly we make sure that the effective interaction is the same across all dimensions, and we should obtain similar equilibrium behavior.

\begin{figure}[h!]
        \centering
         \begin{subfigure}[b]{0.32\textwidth}
             \centering
             \includegraphics[width=\textwidth]{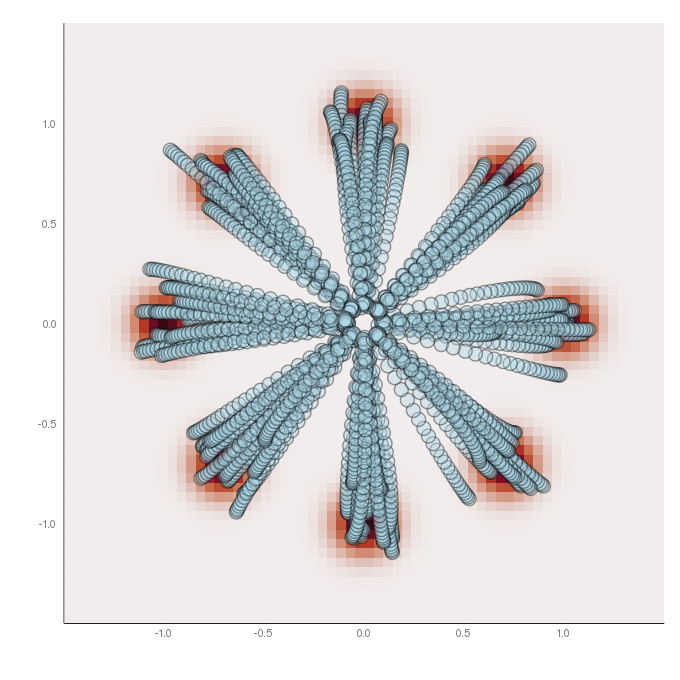}
             \caption{$d=50, \mu=10, \hat\sigma=0.2$.}
             \label{fig:8gaussian-all-d-50-sigma-02}
         \end{subfigure}
         \begin{subfigure}[b]{0.32\textwidth}
             \centering
             \includegraphics[width=\textwidth]{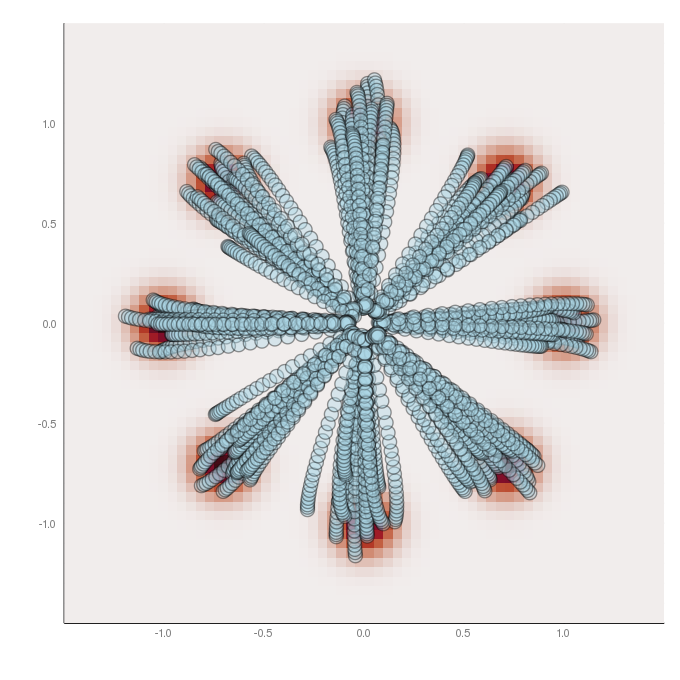}
             \caption{$d=100, \mu=10, \hat\sigma=0.2$.}
             \label{fig:8gaussian-all-d-100-sigma-02}
         \end{subfigure}
        \caption{Agents' trajectories in experiments C for $\hat{\sigma}=0.2$ plotted on the first two dimensions. Agents move from 8 Gaussian distributions (colored red) to the target point ($0,0$). Each plot shows the trajectories solved in different dimensions: (a) $d=50$, (b) $d=100$.}
        \label{fig:8gaussian-all-d-sigma-0.2}
    \end{figure}

\begin{figure}[h!]
        \centering
         \begin{subfigure}[b]{0.32\textwidth}
             \centering
             \includegraphics[width=\textwidth]{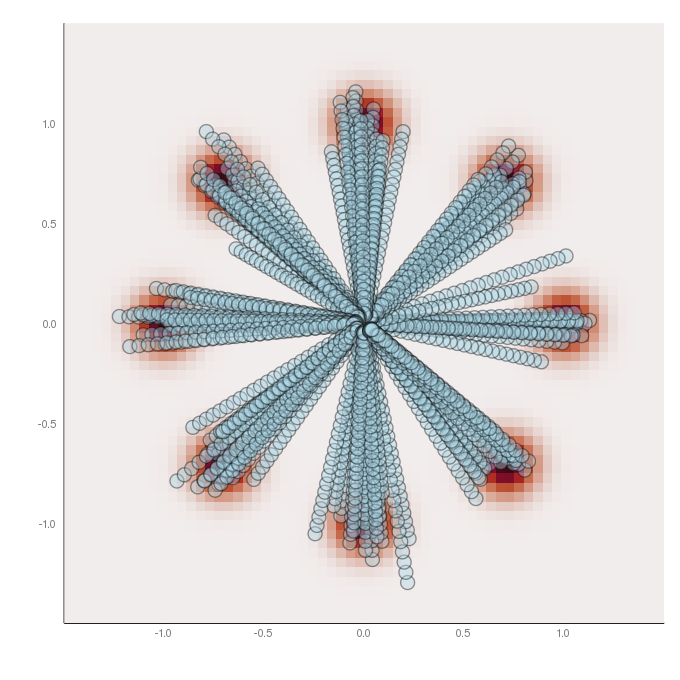}
             \caption{$d=50, \mu=1, \hat\sigma=1.25$.}
             \label{fig:8gaussian-all-d-50-sigma-125}
         \end{subfigure}
         \begin{subfigure}[b]{0.32\textwidth}
             \centering
             \includegraphics[width=\textwidth]{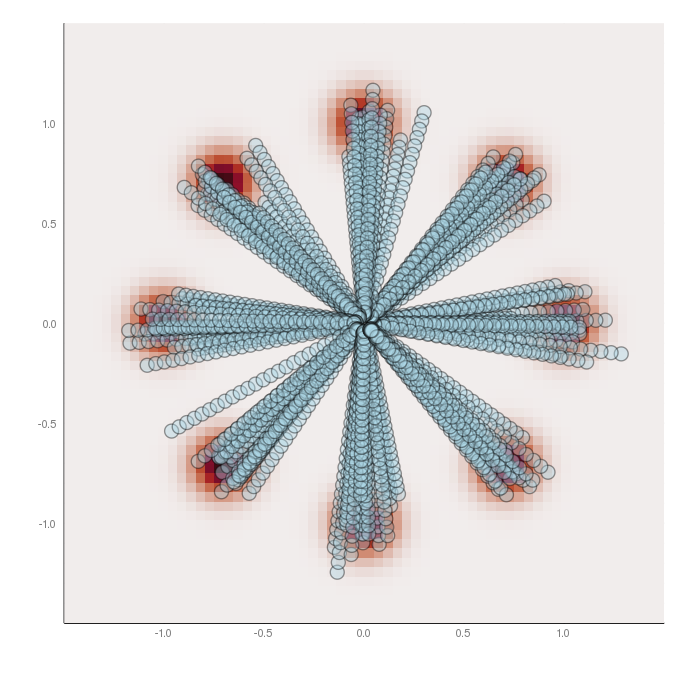}
             \caption{$d=100, \mu=1, \hat \sigma=1.25$.}
             \label{fig:8gaussian-all-d-100-sigma-125}
         \end{subfigure}
        \caption{Agents' trajectories in experiments C for $\hat{\sigma}=1.25$ plotted on the first two dimensions. Agents move from 8 Gaussian distributions (colored red) to the target point ($0,0$). Each plot shows the trajectories solved in different dimensions: (a) $d=50$, (b) $d=100$.}
        \label{fig:8gaussian-all-d-sigma-1.25}
    \end{figure}

The results for $\hat{\sigma}=0.2$ and $\hat{\sigma}=1.25$ are plotted in Figures \ref{fig:8gaussian-all-d-sigma-0.2} and \ref{fig:8gaussian-all-d-sigma-1.25}, respectively.

\begin{figure}[h!]
        \centering
         \begin{subfigure}[t]{0.49\textwidth}
             \centering
             \includegraphics[width=\textwidth]{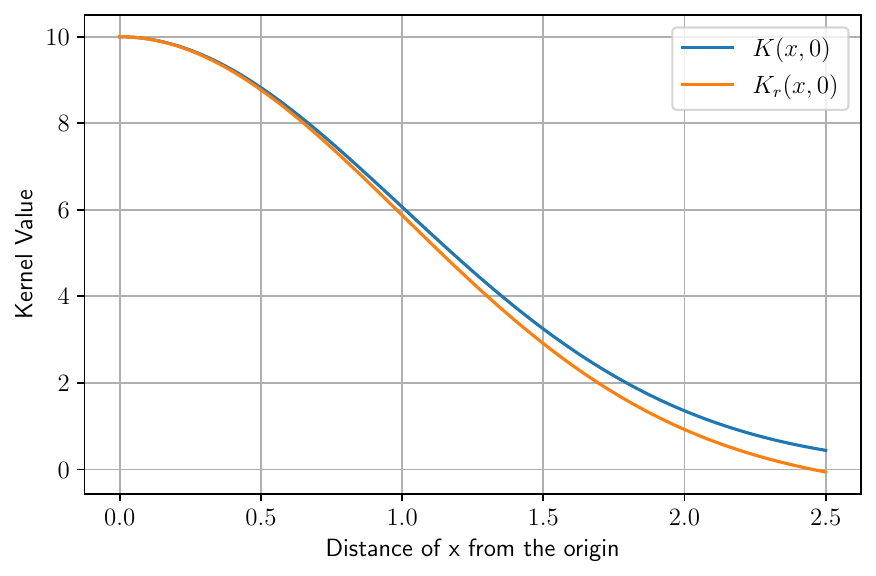}
             \caption{Kernel values along a random direction $r = 512$}
             \label{fig:c1-true-1d-kernel}
         \end{subfigure}
         \begin{subfigure}[t]{0.49\textwidth}
             \centering
             \includegraphics[width=\textwidth]{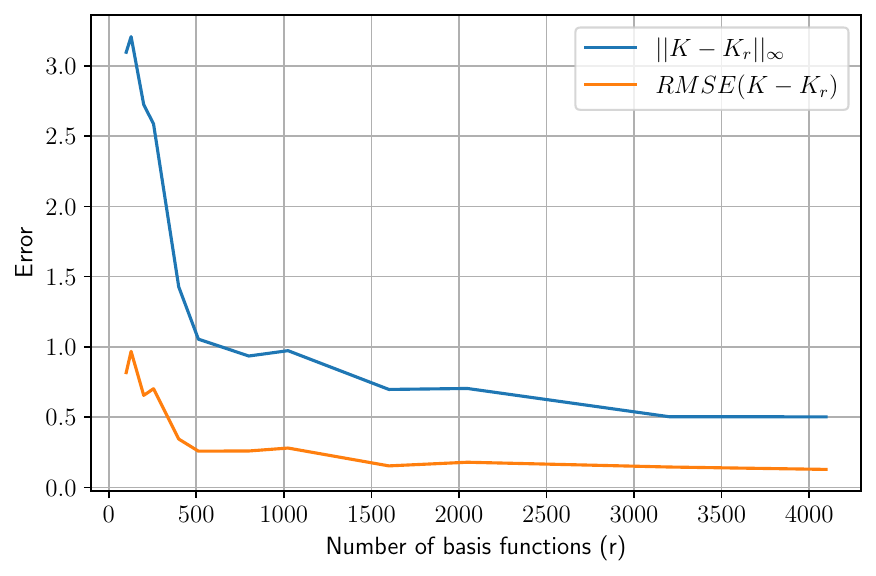}
             \caption{Convergence of errors}
             \label{fig:c1-convergence}
         \end{subfigure} 

        \caption{Kernel approximation for $d=50, \hat{\sigma} = 0.2, \mu = 10.0$.}
        \label{fig:full-interaction-kernel-approx-sigma-02}
    \end{figure}

\begin{figure}[h!]
    \centering
     \begin{subfigure}[t]{0.49\textwidth}
         \centering
         \includegraphics[width=\textwidth]{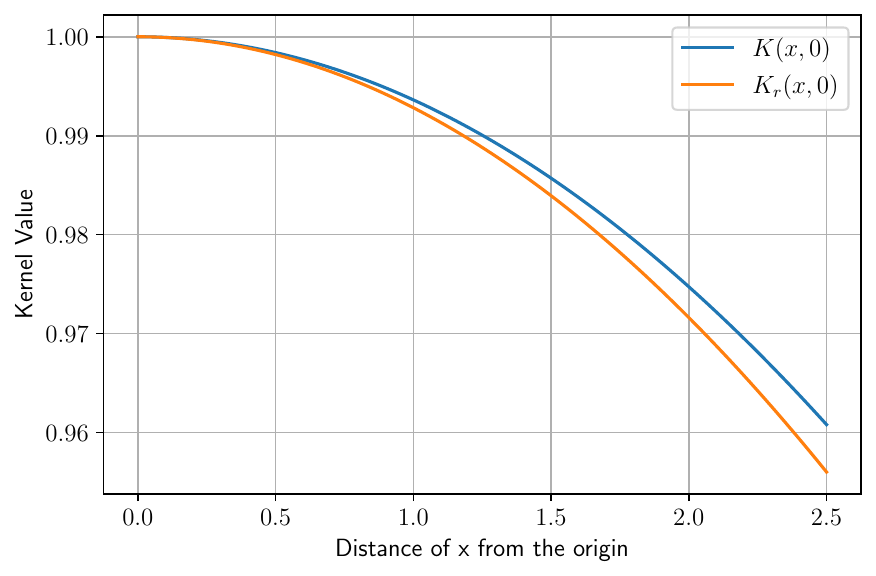}
         \caption{Kernel values along a random direction $r = 512$}
         \label{fig:c2-true-1d-kernel}
     \end{subfigure}
     \begin{subfigure}[t]{0.49\textwidth}
         \centering
         \includegraphics[width=\textwidth]{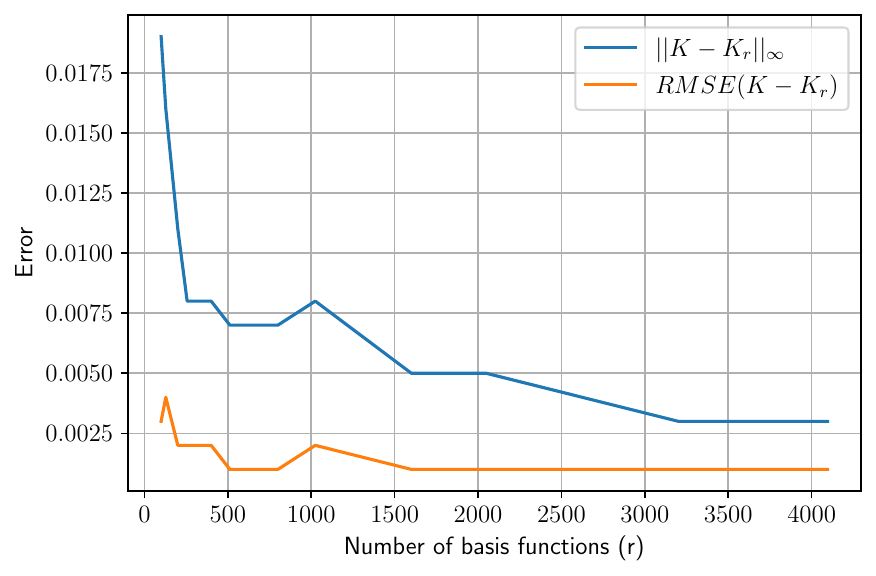}
         \caption{Convergence of errors}
         \label{fig:c2-convergence}
     \end{subfigure} 

    \caption{Kernel approximation for $d=100, \hat{\sigma} = 1.25, \mu = 1.0$.}
    \label{fig:full-interaction-kernel-approx-sigma-1.25}
\end{figure}

\begin{figure}[h!]
    \centering
         \includegraphics[width=0.70\textwidth]{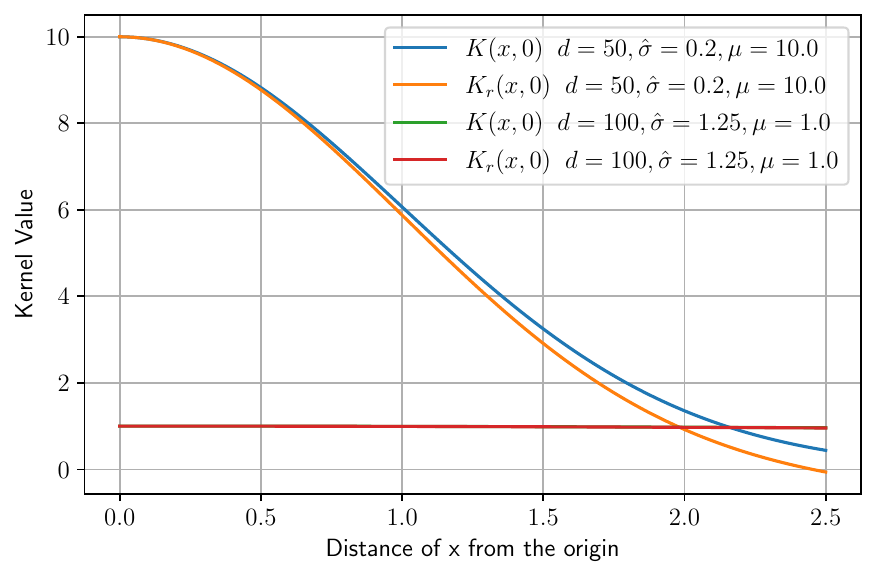}
         \caption{Kernel values along a random ray for $r = 512, \hat{\sigma}=0.2,1.25$}
    \label{fig:c1-c2-1d-kernel}
\end{figure}

Note that the trajectories are almost straight lines when $\hat{\sigma}=1.25$. In Figures \ref{fig:full-interaction-kernel-approx-sigma-02}, \ref{fig:full-interaction-kernel-approx-sigma-1.25}, \ref{fig:c1-c2-1d-kernel} we plot the original and approximate kernels to explain this phenomenon. More specifically, in Figures \ref{fig:c1-true-1d-kernel} and \ref{fig:c2-true-1d-kernel} we plot $K(x,0)$ and $K_r(x,0)$ along a random direction so that $|x| \leq 2.5$. Furthermore, in Figures \ref{fig:c1-convergence}  and \ref{fig:c2-convergence} we plot the decay of the approximation error $K(x,\mu_c)-K_r(x,\mu_c)$ in $l^\infty$ and $l^2$ norms for $x$ sampled according to $\rho_0$, where $\mu_c$ is the center of one of the eight constituent Gaussians of $\rho_0$. Finally, we superimpose Figures \ref{fig:c1-true-1d-kernel} and \ref{fig:c2-true-1d-kernel} in Figure \ref{fig:c1-c2-1d-kernel}.

As we can see in Figures \ref{fig:c2-true-1d-kernel} and \ref{fig:c1-c2-1d-kernel}, the interaction kernel is almost flat for $\hat{\sigma}=1.25$ within the support of $\rho_0$. Hence, the interaction cost is approximately the same for all agents, which effectively decouples the agents leading to individual control problems with a purely quadratic cost. In the latter case, optimal trajectories are straight lines as follows from the Hopf-Lax theory~\cite[Section 3.3]{evansPDE}.

\begin{figure}[h!]
    \centering
     \begin{subfigure}[b]{0.48\textwidth}
         \centering
         \includegraphics[width=\textwidth]{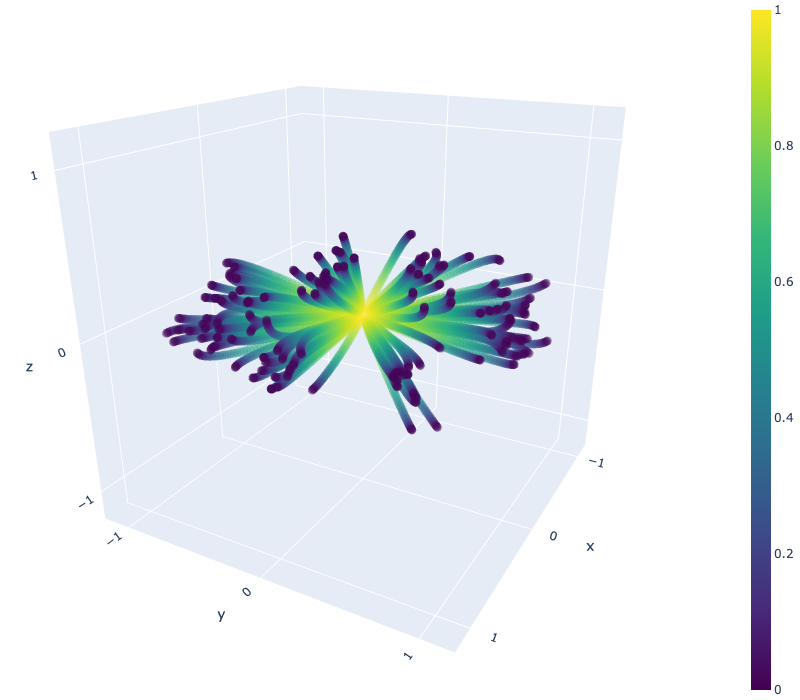}
         \caption{3D view.}
     \end{subfigure}
     \begin{subfigure}[b]{0.48\textwidth}
         \centering
         \includegraphics[width=\textwidth]{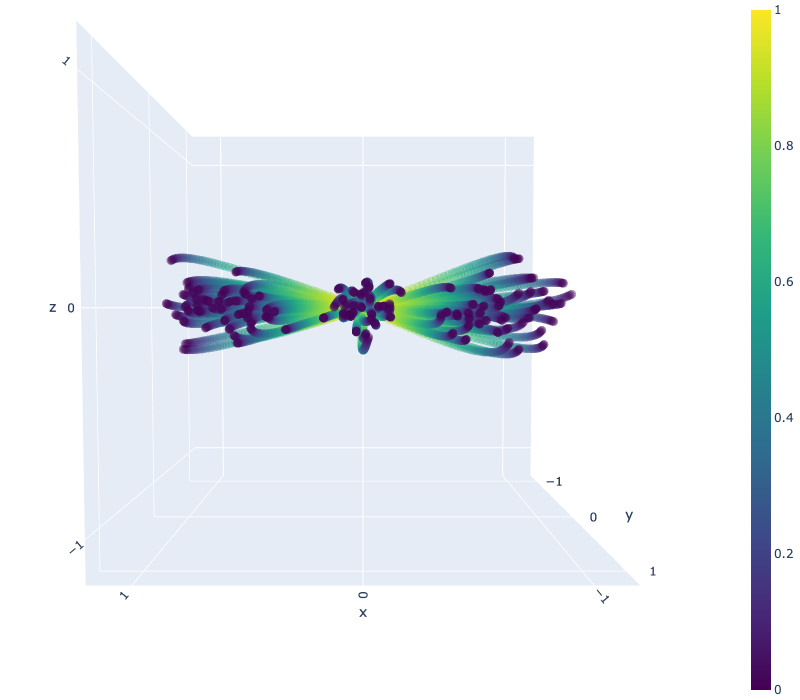}
         \caption{XZ view.}
     \end{subfigure}
     \begin{subfigure}[b]{0.48\textwidth}
         \centering
         \includegraphics[width=\textwidth]{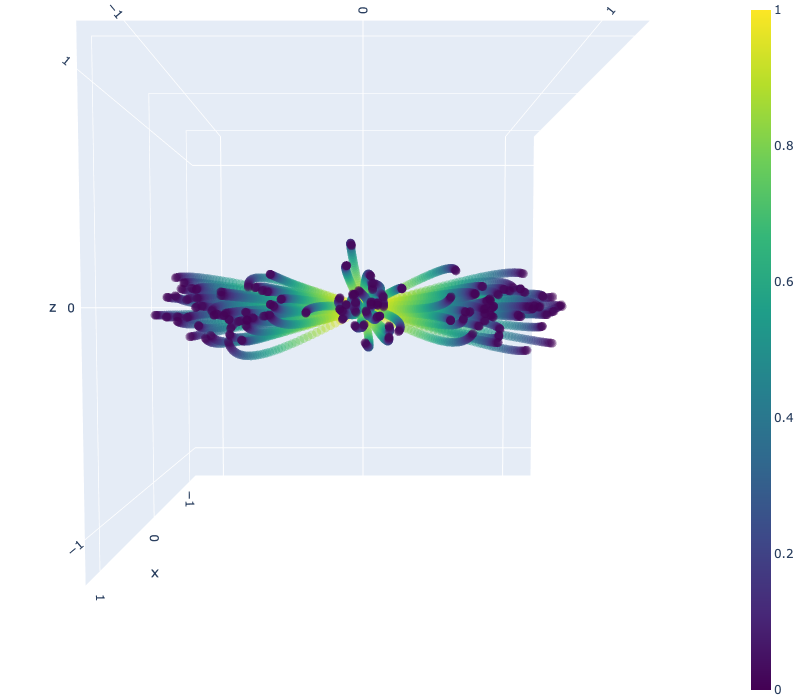}
         \caption{YZ view.}
     \end{subfigure}
     \begin{subfigure}[b]{0.48\textwidth}
         \centering
         \includegraphics[width=\textwidth]{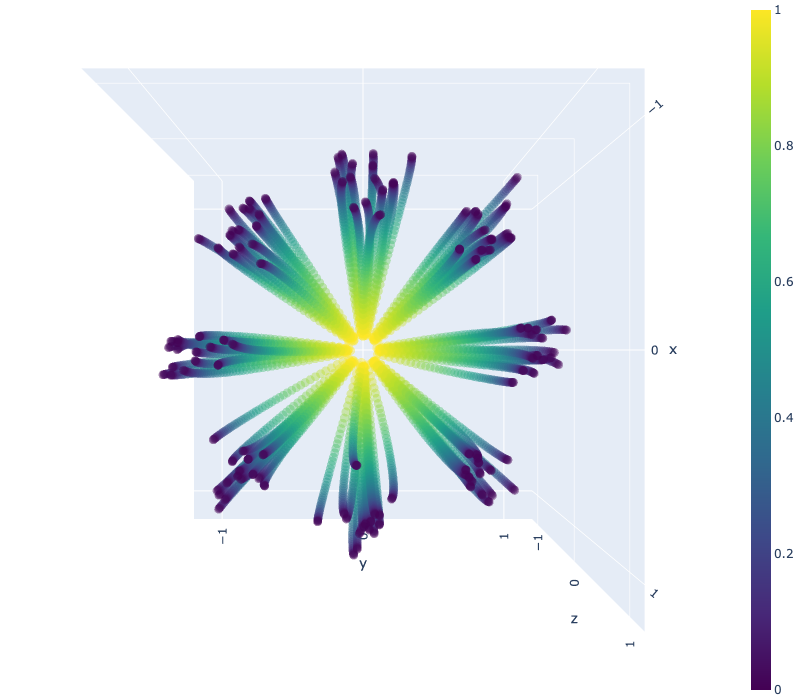}
         \caption{Top view.}
     \end{subfigure}
    \caption{3D plots of Agents' trajectories in experiments C with full-dimensional interactions for $d=50$, $\hat{\sigma}=0.2$. The plots show the first three dimensions of the trajectories. Each agent starts from $t=0$ (colored blue) to $t=1$ (colored yellow). The plots are from four different viewpoints. (a): 3D view, (b), (c): side views (XZ view and YZ view), (d): top view (XY view). Because of full-dimensional interactions, the spread of agents' trajectories can be observed from every viewpoints.}
    \label{fig:8gaussian-all-d-3d}
\end{figure}

We plot the 3D trajectories in Figure \ref{fig:8gaussian-all-d-3d} and report running, interaction, terminal, and total costs in Table \ref{table:C}.

\begin{table}[h!]
    \centering
    \begin{tabular}{|c|c|c|c|c|c|c|}
            \hline
            $d$   & $\hat\sigma$ & $\mu$ & running & interaction & terminal & total
            \\
            \hline
            50 & 0.2  & 10 & 1.05    & 1.96 & 0.0192 &  3.20
            \\
            50 & 1.25  & 1 & 0.674    & 0.492 & 0.00340 &  1.20
            \\
            \hline
            100 & 0.2   &  10 &1.21 & 2.59 & 0.0177 & 3.97
            \\
            100 & 1.25   &  1 & 0.912 & 0.495 & 0.00458 & 1.45
            \\
            \hline
        \end{tabular}
    \caption{Running, interaction, terminal, and total costs in experiments C.}
    \label{table:C}
\end{table}

\section{Conclusion}

We propose an efficient solution approach for high-dimensional nonlocal MFG systems utilizing random-feature expansions of interaction kernels. We thus bypass the costly state space discretizations of interaction terms and allow for straightforward extensions of virtually any single-agent trajectory optimization algorithm to the mean-field setting. As an example, we extend the direct transcription approach in optimal control to the mean-field setting. Our numerical results demonstrate the efficiency of our method by solving MFG problems in up to a hundred-dimensional state space. To the best of our knowledge, this is the first instance of solving such high-dimensional problems with non-deep-learning techniques.

Future work involves the extension of our method to affine controls arising in, e.g., quadrotors~\cite{carrillo2013modeling}, as well as alternative trajectory generation methods that involve deep learning~\cite{ruthotto2020machine,lin2021alternating}.  

Compact feature space representations of interaction kernels are also valuable for inverse problems. In a forthcoming paper~\cite{chow2022numerical}, we recover the interaction kernel from data by postulating its feature space expansion. 

Furthermore, note that feature space expansions of the kernel are not related to the mean-field idealization. Thus, we plan to investigate applications of our method to possibly heterogeneous multi-agent problems where the number of agents is not large enough for the mean-field approximation to be valid~\cite{onken2021multi,onken2021neural}.

Finally, an interesting and challenging question is the convergence analysis of the primal-dual algorithm \ref{eq:algo} described in Section \ref{sec:algorithm}. We anticipate analysis methods developed in~\cite{hadi17b} to be useful for this question.

\section*{Acknowledgments}
\noindent Wonjun Lee, Levon Nurbekyan, and Samy Wu Fung were partially funded by AFOSR MURI FA9550-18-502, ONR N00014-18-1-2527, N00014-18-20-1-2093, and N00014-20-1-2787.

\section*{Appendix}

\begin{proof}[Derivation of \ref{eq:a_inclusion}]
Assume that $v(t,x)$ is a smooth vector field. For every $x\in \R^d$ denote by $z_x(t)$ the solution of the ODE
\begin{equation}\label{eq:flow}
    \dot{z}_x(t)=v(t,z_x(t)),\quad z_x(0)=x.
\end{equation}
If agents are distributed according to $\rho_0$ at time $t=0$ and follow the flow in \ref{eq:flow}, their distribution, $\rho(t,x)$, satisfies the continuity equation
\begin{equation*}
    \partial_t \rho(t,x)+\nabla \cdot \left(\rho(t,x)v(t,x)\right)=0,\quad \rho(0,x)=\rho_0(x).
\end{equation*}
Now assume that $\phi_a$ is the solution of \ref{eq:phi_a}. From the optimal control theory we have that
\begin{equation*}
\begin{split}
    &\int_{\R^d} \phi_a(0,x) d\rho_0(x) \\
    \leq& \int_{\R^d} \left[ \int_0^T \left\{L(t,z_x(t),v(t,z_x(t)))+\sum_{i=1}^r a_i(t)\zeta_i(z_x(t))\right\}dt+\psi(z_x(T))\right] d\rho_0(x) \\
    =&\int_{\R^d}\int_0^T \left\{L(t,x,v(t,x))+\sum_{i=1}^r a_i(t)\zeta_i(x) \right\}  dt d\rho(t,x)+ \int_{\R^d}\psi(x)d\rho(T,x),
\end{split}
\end{equation*}
where equality holds for $(\rho,v)=(\rho_a,v_a)$ given by
\begin{equation}\label{eq:v_rho_a_optim}
\begin{split}
    &v_a(t,x)=-\nabla_p H(t,x,\nabla \phi_a(t,x))\\
    &\partial_t \rho_a(t,x)-\nabla \cdot \left(\rho_a(t,x)\nabla_p H(t,x,\nabla \phi_a(t,x))\right)=0,\quad \rho_a(0,x)=\rho_0(x).
\end{split}
\end{equation}
Summarizing, we obtain that
\begin{equation}\label{eq:phi_a_val_function}
    \begin{split}
        \int_{\R^d} \phi_a(0,x) d\rho_0(x)=&\inf_{\substack{\partial_t \rho+\nabla \cdot (\rho v)=0\\ \rho(0,x)=\rho_0(x)}} \int_{\R^d}\int_0^T L(t,x,v(t,x))  dt d \rho(t,x)\\
        &+\int_{\R^d} \int_0^T \sum_{i=1}^r a_i(t)\zeta_i(x)  \rho(t,x) dt dx + \int_{\R^d}\psi(x)d\rho(T,x),
    \end{split}
\end{equation}
where the equality holds for $(\rho,v)=(\rho_a,v_a)$ in \ref{eq:v_rho_a_optim}. Applying perturbation analysis for optimization problems \cite[Proposition 4.12]{bonnans00} we obtain
\begin{equation}\label{eq:grad_a_appendix}
    \frac{\delta}{\delta a_i(t)}\int_{R^d} \phi_a(0,x)d\rho_0(x)=\int_{\R^d} \zeta_i(x)d\rho_a(t,x)=\int_{\R^d} \zeta_i(z_{x,a}(t))d\rho_0(x),
\end{equation}
where $z_{x,a}$ is the solution of \ref{eq:flow} for the optimal control $v=v_a$.

Now we are in the position for proving the equivalence between \ref{eq:mfg_approx} and \ref{eq:a_inclusion}. We have that
\begin{equation}\label{eq:convolution}
    \int_{\R^d} K_r(x,y)d\rho(t,y)= \sum_{i=1}^r a_i(t) \zeta_i(x),
\end{equation}
where
\begin{equation}\label{eq:a-s}
    a_i(t)=\sum_{j=1}^r k_{ij} \int_{\R^d} \zeta_j(y) d \rho(t,y).
\end{equation}
Therefore, $(\phi,\rho)$ solve \ref{eq:mfg_approx} if and only if $(\phi,\rho)=(\phi_a,\rho_a)$ for $a$ satisfying \ref{eq:a-s}. Furthermore, \ref{eq:grad_a_appendix} yields that \ref{eq:a-s} is precisely equivalent to
\begin{equation*}
    a(t)=\rmK \frac{\delta}{\delta a(t)}\int_{R^d} \phi_a(0,x)d\rho_0(x)
\end{equation*}
which leads to \ref{eq:a_inclusion}.
\end{proof}

\bibliography{literature}

\end{document}